\documentclass{article}
\usepackage{diagrams}

\usepackage{concmath}

\usepackage[T1]{fontenc}
\usepackage{mathptmx}

\usepackage[textwidth=12.8cm, textheight=18.5cm]{geometry}

\usepackage{amssymb}
\usepackage{amsmath}

\newtheorem{prop}{Proposition}
\newtheorem{thm}{Theorem}

\title{Monads and extensive quantities}
\date{}
\author{Anders Kock\\
University of Aarhus}

\newcommand{\E}{\mathcal{E}}
\newcommand{\p}{\pitchfork}
\newcommand{\tto}{\rightrightarrows}
\newcommand {\imes}{\times}
\begin{document}
\maketitle
\small \noindent {\bf Abstract.} If $T$ is a commutative monad on a cartesian 
closed category, then there exists a natural $T$-bilinear pairing 
$T(X)\times T(1)^{X}\to T(1)$ (``integration''), as well as a natural 
$T$-bilinear action $T(X)\times T(1)^{X} \to T(X)$. These data 
together make the  endofunctors $T$ and $T(1)^{(-)}$ (co- and 
contravariant, respectively) into 
a system of extensive/intensive quantities, in the sense of Lawvere. 
A natural monad map from $T$ to a certain monad of distributions (in 
the sense of functional analysis (Schwartz))  arises from the integration.

\normalsize

\section*{Introduction}

Another word for ``extensive quantitiy'', and one which is commonly used outside 
mathematics, is ``distribution" \footnote{In some Germanic 
languages, like German or Danish, the commonly non-mathematical word 
for these kind of distributions is ``Verteilung,'' resp.\  
``fordeling''. In 
mathematics, the word ``distribution'' has acquired a more specific meaning, namely 
the distributions in the sense of Laurent Schwartz, where a 
``distribution'' is a continuous linear 
functional on a space of ``test functions''. }.
In this common non-mathematical usage,  an extensive quantity (say, of rain on a 
given day) 
may be {\em distributed} over a given space, (say the sidewalk), and its 
{\em total} over that space is measured in terms of some absolute quantity: the total 
{\em mass} of the rain on the sidewalk, or the total {\em number} of 
raindrops on the sidewalk (this number is an integer). So we have different quantity {\em 
types} for such totals, say the quantity type ``mass'', or the 
quantity type ``(integral) number''. In this case, both are 
``positive''; but one has also quantity types like ``electric 
charge'' whose quantities admit opposite signs,  in the sense that two such quantities can cancel 
each other. Note that a mass is not a (non-negative real) {\em number}, but 
only becomes so after choosing a {\em unit} of mass. The {\em ratio} of a 
given  mass distribution  with a chosen unit is a  
``(distributed) dimensionless quantity'', and a good 
approximative mathematical model for many types of totals of such 
distributed dimensionless quantities 
is the ring of real 
numbers -- although for mass, say, {\em non-negative} real numbers would be 
a more realistic model for such  quantity.  
Or, the dimensionless quantity may be an integer (or a non-negative integer), 
for the case of ``number of raindrops''. 

A simple  approximative mathematical model of these types of 
dimensionless total quantities 
is: they form commutative (additively written) monoids, like 
${\mathbb R}$, ${\mathbb R}_{+}$, or ${\mathbb N}$, in fact, are 
{\em free} ``algebras''  on one generator, for a suitable notion of 
``algebra'' (thus ${\mathbb N}$ is the free {\em commutative monoid} on one 
generator, and ${\mathbb R}$ is the free {\em real vector space} on one 
generator). The notion of ``free algebra'' may be encoded by the 
notion of monad $T$. Thus  ${\mathbb R}$ is $T(1)$, for the 
free real vector space monad $T$. The fact that $T(1)$ is endowed with a 
$T$-bilinear multiplication comes out from the strength of the monad.

In the present article, we experiment with the viewpoint that the 
dimensionless distributions on a space $X$ form themselves a space 
$T(X)$, where $T$ is a  monad (assumed ``commutative'') on ``the'' cate\-gory $\E$ of spaces, 
assumed to be a cartesian closed category. For instance, $T$ may be 
the ``free commutative monoid'' monad, or the ``free real vector 
space monad'', -- assuming that the field of reals is itself suitably an 
object of $\E$. In fact, we have in mind the case where $\E$ is for 
instance the category of convenient vector spaces and the smooth maps 
in between; or a topos, like the ``smooth topos'', or a ``well 
adapted model'' for synthetic differential geometry; in these 
cases, the cohesion (say, topology) of ${\mathbb R}$ is retained, by ${\mathbb R}$ 
being seen as an object 
of the category $\E$. 

 The units of the monad, i.e.\ the maps $\eta 
_{X}:X\to T(X)$, assign to 
$x\in X$ the distribution with total $1$, 
and concentrated in $x$, in some contexts:  the ``Dirac distribution 
at $x$''.

We shall also have a fragment of a theory of how quantities with a physical 
dimension, like mass, which are 
not pure quantities, fit into the picture. They are likewise 
covariant endofunctors $M$ on $\E$, but $M(1)$, unlike $T(1)$, does 
not carry a natural multiplication or unit; $M$ is, in some sense, 
a ``torsor'' over the appropriate dimensionless-quantity monad $T$.
(In our \cite{MSPQ}, we considered similar torsor structure, but only 
for ``total'' quantities, i.e.\ {\em not} distributed over an 
extended space.)

In most of the present article, we consider only dimensionless quantities.

The theory presented here  implies an attempt to comment on
Schwartz' dictum 
{\em ``les distribution math\'{e}matiques constituent une 
d\'{e}finition math\'{e}matique correcte des distributions 
rencontr\'{e}es en physique.''}
(\cite{Schwartz} p.\ 84) -- but now with ``distribution'' in the 
sense given by  general commutative monads.

A main thing is that $T$ is a {\em covariant 
functor}. An element $P\in T(X)$ is a distribution on 
$X$. We have a unique map $X\to 1$; applying $T(X\to 1)$ to $P$ yields an element in  $T(1)$, the 
{\em total} of the distribution $P$. This covariant feature of 
extensive quantities was 
stressed by Lawvere; in particular, he stressed that distributions 
(in the sense of functional 
analysis) 
are not to be viewed as ``generalized functions'' (functions behave 
contravariantly; they are intensive quantities), but rather are extensive 
quantities, behaving covariantly (at least when restricted to 
distributions of compact support). We shall return to some of 
Lawvere's more specific theory of extensive quantities in the last 
sections.

\medskip
 One main aspect of the theory to be presented 
here is that there is a canonical comparison transformation $\tau$ from  
 the  monad $T$ into a Schwartz 
(double-dualization) type monad $S$ 
associated to $T$. To distinguish, we call the elements $P$ of $T(X)$ 
{\em concrete} distributions, to distinguish them from distributions 
in the sense of functional analysis.    

\medskip

The technical  underpinning of the present 
theory is the theory of strong (= $\E$-enriched) monads on a symmetric monoidal closed 
category  $\E$, developed by the author in a series of articles in the early 1970s, 
\cite{MSMCC}, \cite{DD},  \cite{CCGBCM}, \cite{BCCM}, and \cite{SFMM}. We begin by recalling 
and expanding some of the aspects of this theory; however, since we 
shall be interested in the case of a CCC (= cartesian closed category) $\E$ only, we 
use notation etc. from this special case throughout (so we write 
$\times$ rather than $\otimes$), even though the material in 
Sections 1-3 deal with the general SMC (= symmetric monoidal closed) 
case.

Following the typographical convention from these papers, we write $X\p Y$ for the 
exponential object $Y^{X}$. The counit for the adjunction defining 
the exponential $-\p Y$ is denoted $ev:  X\times (X\p Y ) \to Y$ ( for 
``evaluation'').

Because an unspecified endofunctor $T$ is involved throughout, we 
have in the main preferred to formulate constructions etc.\ in terms of diagrams, 
rather than in terms of elements of the objects considered; however, 
expressions talking about ``elements'' in ``sets'' are sometimes more 
readable than diagrams,  so we shall sometimes use such expressions, even though objects of $\E$ may 
have no (global) elements at all. Alternatively, the elements talked 
about are {\em generalized} elements in the sense used in, say, 
Synthetic Differential Geometry, as expounded in \cite{SDG}. It
 is in principle routine to translate equations and 
constructions, expressed in terms of elements, into commutative 
diagrams.

\medskip
\noindent{\bf Acknowledgements.} The dialectics between extensive and 
intensive quantities, as covariant and contravariant, I learned from 
Lawvere, and this was a leading guideline in the present research. This 
was further spurred by reading Cram\'{e}r's introductory text 
\cite{Cramer} on ``calculus of probabilites,'' which explicitly 
stresses the analogy between probability distributions and mass 
distributions - both important cases of extensive quantities. 

-- I  want to acknowledge several fruitful conversations with Michael 
Wright on these topics.
The diagrams of the article were produced using Paul Taylor's 
``diagrams'' package.

\section{Combinators for strong endofunctors and monads}\label{BCX}

We consider a CCC $\E$; notions of ``enrichment'' or ``strength'' 
refer to this $\E$. 

First, we have the 
evaluation map
$$ev_{X,Y}: X\times (X\p Y)\to Y$$
and its twin sister $\tilde{ev}_{X,Y} :(X\p Y)\times X \to Y$; they are the 
counits for the adjunction $(X\times - )\dashv (X\p -)$ (resp.\ $(-\times 
X) \dashv (X\p -)$).
Often the decorations $X,Y$ on $ev _{X,Y}$ may be omitted, because 
$X$ and $Y$ are clear from the context.

We consider an endofunctor $T:\E \to \E$, assumed strong ( = 
$\E$-enriched); recall 
 that such enrichment is presented in terms of data
$$st_{X,Y}: X\p Y \to T(X)\p T(Y),$$
cf.\ \cite{EK}, or \cite{Borceux} II.6.2.3 item (2).
 
In \cite{MSMCC} and \cite{SFMM}, we observed that the strength can be encoded as a {\em tensorial} strength $$t' _{X,Y}: T(X) 
\times Y \to T(X\times Y),$$ natural in $X$ and $Y$. By 
``conjugating'' with the twist map $X\times Y \to Y\times X$, one 
also gets its ``twin sister''  $$t''_{X.Y}:X\times T(Y) \to T(X\times Y),$$
likewise encoding the strength.

Finally, the strength can be encoded as  {\em cotensorial strength} 
$$\lambda _{X,Y}:T(X\pitchfork Y)\to X\pitchfork T(Y),$$ cf.\ \cite{MSMCC} 
and \cite{CCGBCM}.

We give elementwise descriptions of these basic combinators.
The elements of $X\p Y$ are maps $f:X\to Y$.
The (Eilenberg-Kelly-) strength $st_{X,Y}:X\p Y \to T(X)\p T(Y)$ of $T$ takes such $f$ to $T(f)\in T(X)\p 
T(Y)$. The tensorial strength $t'_{X,Y}:T(X)\times Y \to T(X\times 
Y)$ takes $(P,y)$ to $T(u_{y})(x)$, where $u_{y}:X\to X\imes Y$ takes 
$x$ to $(x,y)$; similarly $t''_{X,Y}:X\imes T(Y) \to T(X\imes Y)$ 
takes $(x,Q)$ to $T(\tilde{u}_{x})(Q)$ where $\tilde{u}_{x}:Y\to X\times 
Y$ takes $y$ to $(x,y)$. Finally, $\lambda _{X,Y}: T(X\p Y) \to X\p 
T(Y)$ takes $S\in T(X\p Y)$ to the map $x\mapsto T(ev_{x})(S)$, where 
$ev_{x}:X\p Y \to Y$ is evaluation at $x\in X$.

\medskip

The following will be a main actor in the following. Let $B\in \E$ 
and let $\beta :T(B)\to B$ be a map. (We are ultimately interested in 
the case where $T$ is a strong monad, and $\beta$ makes $B$ into a 
$T$-algebra). Then for any $X\in \E$, we have  the 
composite
\begin{equation}\label{intx}\begin{diagram}
T(X)\times (X\p B)&\rTo ^{ t'}&T(X\times (X\p B))&\rTo 
^{T(ev)}&T(B)&\rTo ^{\beta }&B.
\end{diagram}\end{equation}
Alternatively, by (\ref{deltay}) below, this map equals
\begin{equation}
\label{intbxx}\begin{diagram}T(X)\times (X\p B) &\rTo ^{id \times 
st}& T(X)\imes (T(X)\p T(B))&\rTo ^{ev}&T(B)&\rTo ^{\beta}&B.
\end{diagram}\end{equation}
In elementwise terms: If $P\in 
T(X)$ and  $\phi \in X\p B$ (so $\phi :X \to B$ is a map), the value 
of (\ref{intbxx}) on $(P,\phi )$ is $\beta (T(\phi )(P)) \in B$, and is denoted $\int 
_{X,B}\phi (x)\; dP(x)$ (with `$x$' as a dummy variable). 
Frequently,  $B$ and $\beta$ may be understood from the context (and the most 
important case is when $B=T(1)$), in which case the definiton of 
$\int$ reads 
\begin{equation}\label{intelx}
\int _{X}\phi (x) \; dP(x) : = \beta (T(\phi 
)(P)),\end{equation}
or, with increasing pedantry (rarely needed)
$$\int _{X}\phi (x) \; dP(x)=\int _{X,B}\phi (x) \; dP(x)=\int 
_{X,(B,\beta )}\phi (x) \; dP(x).$$

We are interested in the case where $\beta$ provides $B$ with 
structure of $T$-algebra. (In this case, when $\E$ is the category of 
sets, this ``integration'' relationship between monad/algebra theory, and  
algebraic theories, dates back to the early days of monad theory with Linton, 
Wraith and others, in the mid 1960s: they knew that the elements of 
$T(X)$ can be interpreted as 
$X$-ary operations $X\p B \to B$ on $T$-algebras $B$.) 
If in particular $X$ is $B$ itself, and we put $\phi = id_{B}$,
we end up with 
\begin{equation}\int _{B,B}x\; dP(x) = \beta (P). 
\label{intelxx}\end{equation}
We are ultimately to read this as the {\em expectation} of $P$, see 
Section \ref{EOMX} below.

We collect some further definitions and 
basic relations concerning the combinators related to a strong endofunctor 
$T:\E \to \E$ on a symmetric monoidal closed category $\E$. We 
continue to use 
notation as if $\E$ were actually cartesian closed, i.e.\ we write 
$\times$ rather than $\otimes$.

 We already considered the counit $ev$  for the 
adjunction $( X\times - )\dashv (X\p -)$. The unit for this 
adjunction is not used so 
often, it is denoted $u$, with suitable decorations. Similarly,  $\tilde{u}$ is the unit 
corresponding to the counit $\tilde{ev}$. Thus 
$$u_{X,Y}:Y\to X\p (X\imes Y) \quad \tilde{u}_{X,Y}: Y\to X\p 
(Y\times X).$$
 The decorations are usually omitted from notation; even the tilde may often 
safely be omitted; one 
case where it is useful to retain the tilde is in the characterizing diagram 
for $\delta$ (canonical map to double dual; ``$\delta$'' for ``Dirac''); we have that the 
following diagram commutes
\begin{equation}\label{deltax}\begin{diagram}
X\times (X\p B)&\rTo ^{\delta \times id}&((X\p B)\p B)\times (X\p B)\\
\dTo^{ev}&&\dTo_{\tilde{ev}}\\
B&\rTo _{=}&B. 
\end{diagram}\end{equation}
Next, a diagram relating the tensorial strength $t'$ with the 
Eilenberg-Kelly strength (= $\E$-enrichment) $st$ of $T$:
\begin{equation}\label{deltay}\begin{diagram}
TA\times (A\p B)&\rTo ^{t'_{A,A\p B}}& T(A\times (A\p B))\\
\dTo^{TA\times st_{A,B}}&&\dTo_{T(ev)}\\
TA\times (TA\p TB)&\rTo_{ev}& TB.
\end{diagram}\end{equation}
The proof of this comes about, via manipulation by exponential 
adjointness, of the definition of $t'_{X,Y}
$ (as given in 
\cite{MSMCC} p.2 for $t''$)
 in terms of $st$, namely as the exponential adjoint of the composite
$$\begin{diagram}Y&\rTo ^{u_{X,Y}}&X\p (X\times Y)&\rTo ^{st}&TX\p 
T(X\times Y).
\end{diagram}$$
(put $Y:=A\p B$ and $X:=A$).

Next, a diagram relating the tensorial strength $t'$ with the 
cotensorial strength $\lambda$: 
\begin{equation}\label{deltaz}\begin{diagram}
T(Y\p Z)\times Y&\rTo ^{\lambda_{Y,Z} \times Y}& (Y\p TZ)\times Y\\
\dTo ^{t'_{Y\p Z, Y}}&& \dTo_{\tilde{ev}}\\
T((Y\p Z)\times Y)&\rTo _{T(\tilde{ev})}& TZ.
\end{diagram}\end{equation}
This also follows by exponential adjointness manipulations; a proof 
is given in \cite{CCGBCM} Lemma 1.2.

We shall also have occasion to use commutativity of the outer diagram 
in
\begin{equation}\label{intxx}\begin{diagram}[nohug]
TX\times (X\p B)&\rTo ^{T \delta \times id}&T((X\p B)\p B)\times (X\p 
B)& \rTo ^{\lambda \times id}& ((X\p B)\p TB)\times (X\p B)\\
 &&\dTo_{t'}&&\\
\dTo^{t'}&&T((X\p B)\p B)\times (X\p B))&&\dTo_{\tilde{ev}}\\
&\ruTo_{T(\delta \times id)}&&\rdTo_{T(\tilde{ev})} &\\
T(X\times (X\p B))&&\rTo _{T(ev)}&&TB;
\end{diagram}\end{equation}
here, the left hand square commutes by naturality of $t'$, the right 
hand square by (\ref{deltaz}), and the triangle by applying $T$ to 
(\ref{deltax}).

We will take the tensorial strength $t'$ (or equivalently $t''$) as the primary encoding. If 
$\E$ is the category of sets, $t'_{X,Y}:T(X)\times Y \to T(X\times Y)$ is the 
map which for $P\in T(X)$ and $y\in Y$ returns the value $T(u_{y})(P) 
\in T(X\times Y)$, where $u_{y}:X\to X\times Y$ is the map 
$x\mapsto (x,y)$.

The combinator $t'$ satisfies a unit law and an associative law. The 
unit law says that
$t'_{X,1}:T(X)\times 1 \to T(X\times 1)$ is the composite of the two 
canonical isomorphisms  $T(X)\times 1 \cong T(X)\cong T(X\times 1)$. 
The associative law says that the composite
$$\begin{diagram}T(X)\times Y \times Z &\rTo ^{t'_{X,Y}\times 
Z}&T(X\times Y)\times Z&\rTo ^{t'_{X\times Y,Z}}&T(X\times Y\times 
Z)\end{diagram}$$ equals $t'_{X,Y\times Z}$. (For simplicity, we 
write as if $\times$ were strictly associative.) There are  similar 
unit- and associative laws for $t''$. All these laws follow from the 
standard laws for the $\E$-enrichment, cf.\ \cite{MSMCC}.

We shall have occasion to use a ``derived'' combinator,
\begin{equation}\label{mtx}t_{X,Y,Z}: X\times T(Y) \times Z \to 
T(X\times Y \times Z);\end{equation}
it can be defined in several ways which are equivalent in view of the 
associative law for $t'$ and the construction of $t''$ in terms of 
$t'$. One way to define it is to consider
\begin{equation}\label{txx}t''_{X\times Z,Y}:X\times Z \times T(Y) 
\to T(X\times Z\times Y)\end{equation}
and conjugate it by interchange of  $Z$ and $T(Y)$, respectively, $Z$ 
and $Y$. It can, by the 
associative law, also be defined as the composite
\begin{equation}
\label{tx}\begin{diagram}
X\times T(Y)\times Z&\rTo ^{X\times t'_{Y,Z}}&X\times T(Y\times 
Z)&\rTo ^{t''_{X,Y\times Z}}&T(X\times Y \times Z).$$
\end{diagram}\end{equation}

A natural transformation $\tau: T\Rightarrow S$ between two 
strong functors is {\em strong} if all squares of the form
\begin{equation}\label{strx1}\begin{diagram}T(X)\times Y&\rTo ^{t'_{X,Y}}&T(X\times Y)\\
\dTo^{\tau _{X}\times Y}&&\dTo_{\tau _{X\times Y}}\\
S(X)\times Y&\rTo_{s'_{X,Y}}&S(X\times Y)
\end{diagram}\end{equation}
commute, where $t'$ and $s'$ are the tensorial strengths of $T$ and 
$S$, respectively\footnote{There is another use of the word 
``strong'' for a natural transformation, namely ``all the naturality 
squares are pull-backs''. This is not how we use the word here.}.

Let $(T,\eta, \mu , t')$ be a strong monad, (so $T$ is equipped with a 
strength $t'$, and $\eta$ and $\mu$ are strong natural 
transformations; recall that $t'$ induces a natural strength on 
$T\circ T$; ``strengths compose''; an explicit expression for the 
composite strength (in $t''$-terms) appears in the center line of (\ref{starx}) below). 

It is easy to deduce  from the strength of the natural transformation 
$\eta :id \Rightarrow T$ that the following diagram commutes:
\begin{equation}
\label{intxy}\begin{diagram}X\imes (X\p B)&&\rTo ^{ev}&&B\\
\dTo^{\eta_{X}\times id}&&&&\dTo_{\eta_{B}}\\
T(X)\times (X\p B)&\rTo_{t'}&T(X\times (X\p B))&\rTo_{T(ev)}&T(B).
\end{diagram}\end{equation}

Together with a monad $T$ on $\E$ comes the notion of 
(Eilenberg-Moore-) $T$-{\em algebra} $(B,\beta )$ where $\beta 
:T(B)\to B$ satisfies a unit- and associative law. In particular 
$\beta \circ \eta _{B} =id_{B}$. So from the above, we deduce that if 
$B=(B,\beta )$ is a $T$-algebra, then
 precomposing $\int _{X}:T(X)\times (X\p B)\to B$ with $\eta _{X}\times id$  
just yields the evaluation map (use the description (\ref{intx})); let us record this:
\begin{equation}\label{intxz}\begin{diagram}\bigl[X\times (X\p B) 
&\rTo^{\eta_{X}\times id}& T(X)\times (X\p B)&\rTo^{\int_{X}}& 
B\bigr] =&\bigl[ X\times (X\p B)&\rTo ^{ev}&B\bigr].
\end{diagram}\end{equation}
The $T$-algebras form a category 
$\E ^{T}$, whose maps 
are called {\em $T$-homomorphisms}; we shall also use the term {\em 
$T$-linear} maps, because this will allow us to talk about 
$T$-{\em bilinear} maps, a notion introduced in the strong-monad 
context in \cite{Linton} and 
\cite{BCCM}. We shall recall and expand some of the theory from loc.cit.
As long as the monad $T$ is fixed,  we may say ``linear'' instead of $T$-linear, 
and similarly for ``bilinear''.

The following result (Theorem 2.1 in 
\cite{MSMCC})
 is important for our present aims: 
\medskip

\noindent{\em the functor part $T$ of a strong monad carries two canonical 
structures as a monoidal functor; with respect to each of these, $\eta$ is a 
monoidal transformation.}

\medskip

These two monoidal structures are in loc.cit.\ denoted $\psi$ and 
$\tilde{\psi}$, respectively; one is just the ``twisted'' version of 
the other; $\psi _{X,Y}$ is the composite
\begin{equation}\label{phix}\begin{diagram}T(X) \times T(Y) & \rTo ^{t'_{X,T(Y)}}& T(X\times TY)&\rTo 
^{T(t''_{X,Y})}&T^{2}(X\times Y)&\rTo ^{\mu _{X\times Y}}&T(X\times Y)
\end{diagram},\end{equation}
and $\tilde{\psi}_{X,Y}$ similarly is
\begin{equation}\label{phixx}\begin{diagram}
T(X)\times T(Y)&\rTo ^{t''_{T(X),Y}}&T(TX\times Y)&\rTo 
^{T(t'_{X,Y})}&T^{2}(X\times Y)&\rTo^{\mu _{X\times Y}}&T(X\times Y).
\end{diagram}
\end{equation}

The ``nullary'' part of both the monoidal structures is $\eta _{1}: 
1\to T(1)$, where $1$ is the unit object of $\E$ (i.e.\ the terminal 
object, in the CCC case).

In distribution theory, if $P$ an $Q$ are distributions of compact 
support on spaces $X$ and $Y$ respectively, the distribution $\psi 
_{X,Y}(P,Q)$ on $X\times Y$ is called the {\em tensor product} of $P$ 
and $Q$, cf.\ \cite{Schwartz} III.1.

\medskip

Recall from \cite{MSMCC} that the strong monad $T$ is called {\em 
commutative} if $\psi = \tilde{\psi}$. 
In Theorem 3.2 in 
\cite{MSMCC} it is proved that if $T$ is commutative in this sense, then 
$\mu$ is a monoidal transformation (and hence $T$ a monoidal monad, 
since $\eta$ is in any case a monoidal transformation).

 (There is a converse result contained in Theorem 2.3 in \cite{SFMM}; 
it contains the assertion that the strength $t'$ of $T$ can be 
reconstructed from the structure of monoidal monad; however, in the 
present article, we prefer to have strength as a principle underlying 
everything -- almost a part of the logic. For $\E =$ Sets, strength 
is canonically present in all functors and transformations. 

\medskip 

We proceed to describe some of the relationships we need between the various 
combinators associated to strong monads.

\begin{prop}Precomposing $\psi _{X,Y}$ with $\eta _{X}\times T(Y)$ 
yields $t''_{X,Y}$. Similarly, precomposing $\tilde{\psi}_{X,Y}$ with 
$T(X) \times \eta _{Y}$ yields $t'_{X,Y}$.
\label{precomp}\end{prop}
Diagrammatically, the first assertion says that the outer diagram in the following 
diagram commutes; the inner square commutes by naturality, the left 
hand triangle commutes since $\eta$ is strong, and the right hand 
triangle commutes by a monad law. So the total diagram is likewise 
commutative, and this proves the first assertion of the Proposition.
$$\begin{diagram}[nohug]&& X\times T(Y)&\rTo ^{t''_{X,Y}}&T(X\times Y)&&\\
&\ldTo ^{\eta _{X}\times TY}&\dTo _{\eta _{X\times TY}}&&\dTo ^{\eta 
_{T(X\times Y)}}&\rdTo ^{id}&\\
TX\times TY&\rTo_{t'_{X,TY}}&T(X\times TY)&\rTo 
_{T(t''_{X,Y})}&T^{2}(X\times Y)&\rTo_{\mu _{X\times Y}}&T(X\times Y).
\end{diagram}$$
The proof and diagram for the second assertion are similar.

\medskip

Let $A=(A,\alpha )$ and $C=(C,\gamma )$ be  $T$-algebras. A map 
$f:A\times X \to C$ is called {\em 1-linear} (or {\em linear} 
($T$-linear) {\em in the first 
variable}), cf.\ \cite{BCCM}, if the following pentagon 
commutes
\begin{equation}\begin{diagram}T(A)\times X &\rTo ^{t'_{A,X}}& T(A\times X)&\rTo 
^{T(f)}&T(C)\\
\dTo ^{\alpha \times X}&&&&\dTo _{\gamma}\\
A\times X&&\rTo _{f}&&C.
\end{diagram}\label{1linx}\end{equation}
Similarly, a map $X\times A \to C$ is called {\em 2-linear} (or 
{\em linear in the second variable}) if a similar diagram, now using $t'' 
_{X,A}: X\times T(A) \to T(X\times A)$ commutes.
Finally, if further $B=(B,\beta)$ is a $T$-algebra, a map $A\times B \to C$ 
is called {\em bilinear} if it is both 1-linear and 2-linear. 

(One may 
define the notion of $n$-linear map $A_{1}\times \ldots A_{n}\to C$ 
(where the $A_{i}$s and $C$ are (underlying objects of) algebras), 
and in this way, one should get a multicategory; however, to 
substantiate this, there are 
some coherence conditions that need to be worked out.)

Recall that an object of the form $T(Z)$ carries a canonical algebra 
structure, namely with structure map $\mu _{Z}:T^{2}(Z)\to T(Z)$ 
(this is the {\em free} $T$-algebra 
on $Z$). The algebras in the following Proposition are free.

\begin{prop}\label{two}For any $X$ and $Y$ in $\E$, the map $t'_{X,Y}:T(X)\times 
Y \to T(X\times Y)$ is 1-linear. Similarly $t'' _{X,Y}:X\times T(Y) 
\to T(X\times Y)$ is 2-linear. 
\end{prop}
{\bf Proof.} The pentagon (\ref{1linx}) above, with $A=T(X)$ and 
$X=Y$ has as top line the map 
$T(t'_{X,Y})\circ t'_{T(X),Y}$, and this is an instance of the tensorial strength 
for the composite functor $T\circ T$; and then the commutativity of the 
pentagon is seen to be an instance of the assumption that $\mu$ is a 
strong natural transformation. The proof of the second assertion 
follows by suitable conjugation by  twist maps.

\medskip

A consequence (cf.\  \cite{BCCM}) is that $\psi _{X,Y}:T(X)\times T(Y) \to 
T(X\times Y)$ is 1-linear, and that $\tilde{\psi}_{X,Y}$ is 2-linear.
If all instances of $\psi$ are 2-linear (or equivalently, bilinear), then the monad is 
commutative, and vice versa, cf.\ loc.cit.\ Proposition 1.5.

\medskip

Recall that if $C=(C,\gamma )$ is a $T$-algebra, any map $X\to C$ 
extends uniquely over $\eta _{X}:X\to T(X)$ to a linear map $T(X)\to C$; this is 
the ``free'' property of $T(X)$. We have a closely related 
``universal'' property of  
 $T(X)\times Y$:
 \begin{prop}\label{twop}Any map $f:X\times Y \to 
C$ extends uniquely over  $\eta _{X}\times Y: X\times Y \to 
T(X)\times Y$ to a 1-linear map : $T(X)\times Y \to C$.\end{prop} For, there are 
natural bijective correspondences
$$\hom (X\times Y,C)\cong\hom (X,Y\p C) \cong \hom_{T}(T(X), Y\p C)$$
(where the second occurrence of $Y\p C$ is the cotensor  $Y\p C $ in $\E ^{T}$, recalled in 
(\ref{cotensor}) below, and where the second bijection is induced by 
precomposition by $\eta _{X}$); finally, the 
set $\hom_{T}(T(X), Y\p C)$ is in bijective correspondence with the 
set of 1-linear maps $T(X)\times Y \to C$, by \cite{BCCM} Proposition 1.3 
(i).

\medskip
It is useful to be have an explicit formula for the 1-linear 
extension of $f:X\times Y \to C$; it is the composite
\begin{equation}\label{extx}
\begin{diagram}
T(X)\times Y & \rTo ^{t'_{X,Y}}&T(X\times Y)&\rTo 
^{T(f)}&T(C)&\rTo ^{\gamma}&C.
\end{diagram}\end{equation}
For, $t'_{X,Y}$ is 1-linear, by Proposition \ref{two}, and the two 
other maps in (\ref{extx}) are linear, so the composite is 1-linear. 
Also, it is easy to see that the restriction of (\ref{extx}) along 
$\eta _{X }\times Y$ gives $f$ back (use the unit law $t'_{X,Y}\circ (\eta _{X}\times 
Y)= \eta _{X\times Y}$, and also $\gamma \circ \eta _{C}= id _{C}$). So 
(\ref{extx}) satisfies the two conditions in Proposition \ref{twop}.

\begin{prop}\label{charpsi}The map $\psi _{X,Y}:T(X)\times T(Y)\to T(X\times Y)$ is 
characterized by the following two properties: it is 1-linear, and its 
precomposition with $\eta _{X}\times TY$ is $t''_{X,Y}$. Similarly 
$\tilde{\psi}:T(X)\times T(Y)\to T(X\times Y)$ is characterized by 
the two properties: it is 2-linear, and its precomposition with 
$TX\times \eta _{Y}$ is $t'_{X,Y}$.
\end{prop}
{\bf Proof.} We prove the first assertion. We already observed that $\psi _{X,Y}$ satisfies these 
two conditions, cf.\ Propositions \ref{precomp} and the quotation 
from \cite{BCCM} after Proposition \ref{two}. The converse follows from 
Proposition \ref{twop}. 

\medskip
Assume now that $B=(B,\beta)$ is an algebra, and consider a map 
$f:X\times B \to C$ which is 2-linear. It extends, by the above, to a 
1-linear map $\overline{f} :T(X)\times B \to C$, and we may ask whether this 
$\overline{f}$ inherits from $f$ the property of being 2-linear (and 
is thus bilinear). A sufficient condition is commutativity of $T$: 

\begin{prop}\label{fivex}Let $T$ be commutative. Let $B=(B,\beta )$ and $C=(C,\gamma )$ be $T$-algebras, 
and assume that $f :X\times B \to C$ is 2-linear. Then its 1-linear 
extension $\overline{f}: T(X)\times  B \to C$ is bilinear.
\end{prop}
{\bf Proof.} We use the formula (\ref{extx}) with $Y=B$ for the 
extension. It remains to prove 2-linearity of this map, i.e.\ to 
prove commutativity of the following diagram (where the bottom line is 
$\overline{f}$, according to (\ref{extx}))
$$\begin{diagram}TX \times TB&\rTo ^{t''}&T(TX\times B)&\rTo 
^{T(t')}&T^{2}(X\times B)&\rTo ^{T^{2}f}& T^{2}C&\rTo^{T\gamma}&TC \\
\dTo^{TX\times \beta }&&&&&&&&\dTo_{\gamma}\\
TX\times B &&\rTo_{t'}&T(X\times B)&\rTo_{Tf}&TC&\rTo_{\gamma }&&C.
\end{diagram}$$
Consider the composite $\gamma \circ 
T\gamma \circ T^{2}f$ of the last three arrows in the clockwise 
composite. By pure monad- and algebra theory, we have  
\begin{equation}\label{replx}\gamma \circ 
T\gamma \circ T^{2}f=\gamma 
\circ Tf \circ \mu ,\end{equation}  and having performed this 
replacement, the definition of  $\tilde{\psi}$ appears at the 
beginning of the clockwise composite. Since $T$ was assumed 
commutative, we may replace $\tilde{\psi}$ by $\psi$, and after this 
replacement, the clockwise composite comes out as the composite
$$\begin{diagram}TX \times TB&\rTo^{t'}&T(X\times TB)&\rTo 
^{T(t'')}&T^{2}(X\imes B)&\rTo ^{\mu}&T(X\imes B)&\rTo ^{Tf}&TC&\rTo 
^{\gamma}&C.
\end{diagram}$$
Now we can use (\ref{replx}) once more, in the opposite direction, 
and we end up with the composite
$$\begin{diagram}TX\times TB&\rTo^{t'}&T(X\times TB)&\rTo 
^{T(t'')}&T^{2}(X\times B)&\rTo ^{T^{2}f}&T^{2}C& \rTo ^{T\gamma}&TC&\rTo^{\gamma}&C .
\end{diagram}$$
After these manipulations with the clockwise composite, the  
diagram to be proved commutative has the following shape
$$\begin{diagram}TX\times TB&\rTo ^{t'}&T(X\imes TB)&\rTo 
^{T(t'')}&T^{2}(X\times B)&\rTo^{T^{2}f}&T^{2}C&\\
\dTo^{TX \times \beta}&&\dTo_{T(X\imes \beta)}&&*&&\dTo _{T\gamma }\\
TX\times B&\rTo_{t'}&T(X\times B)&&\rTo_{Tf}&&TC&&\rTo_{\gamma}&C
\end{diagram}$$
Here the pentagon * commutes: it comes about by applying the functor 
$T$ to the diagram expressing the assumption that $f:X\times B \to C$ is 2-linear.
This proves the desired 2-linearity of $\overline{f}$.

\medskip 

For a commutative $T$, we get as an immediate corollary that a map $f:X\times Y\to C$ 
(where $C=(C,\gamma)$ is a $T$-algebra) 
extends uniquely to a bilinear $T(X)\times T(Y) \to C$. Since also 
$f$ extends uniquely to a linear $T(X\times Y) \to C$, we may 
conclude that $T(X\times Y)$ may serve as $T(X)\otimes T(Y)$ in 
$\E^{T}$, with $\psi _{X,Y}$ as the universal bilinear map;  but we shall not prove or need existence of such tensor 
products for general $T$-algebras.

\medskip
The following is hardly surprising,
and the routine proof is  omitted:
 \begin{prop}\label{montrans} If $\tau $ is a strong natural 
transformation from one strong monad $T$ to another one, $S$, 
compatible with the monad structures, then $\tau$ will also be 
compatible with the monoidal strucures $\psi ^{T}$, $\psi ^{S}$, 
i.e.\ it will be a monoidal transformation. Similarly, $\tau$ will be 
compatible with the monoidal structures $\tilde{\psi}^{T}$, 
$\tilde{\psi}^{S}$.\end{prop}

\section{Monads and double dualization}\label{MDD}
Given a commutative monad $T$ on $\E$. Some of the formal properties of the construction $\int _{X}\phi 
(x)\; dP(x)$ is best stated in terms of a transformation $\tau$ from 
$T$ to a certain ``double dualization'' monad associated to $T$. In 
essence, $\tau _{X}$ will be exponential adjoint of $ \int _{X} : 
T(X)\times (X\p B)\to B$ (where  $B$ is a $T$-algebra).

We  assume that $\E$ has equalizers (or sufficiently many 
-- only a few are needed; we study this question in more detail in 
Section \ref{exax}). In this case, the 
category  $\E ^{T}$ of algebras for a strong monad $T=(T,\eta ,\mu , 
t')$ becomes enriched over $\E$: if $(A,\alpha )$ and $(B ,\beta )$ 
are algebras, the $\E$-valued hom-object $[(A,\alpha), (B; \beta 
)]_{T}$  is carved out of $A\pitchfork B$ by an evident equalizer 
diagram involving $\alpha$ and $\beta$, expressing the $T$-homomorphism condition diagrammatically.
(This goes back to \cite{Bunge}.)
We write $(A, \alpha )\pitchfork _{T}(B,\beta )$ for this hom object, 
and often omit $\alpha $ and  
$\beta $ from notation; they are then to be understood from the 
context.
Note that $(A,\alpha )\pitchfork _{T}(B,\beta )$ is a subobject of 
$A\pitchfork B$. In short notation, $A\p _{T}B \subseteq A\p B$.
 
Also, $\E ^{T}$ is cotensored over $\E$: if $X\in \E$ and $(B,\beta )\in \E ^{T}$, 
the cotensor $X\pitchfork (B, \beta )$ is the object $X\pitchfork B$ 
in $\E$, equipped with the $T$-structure
\begin{equation}\begin{diagram}T(X\pitchfork B)&\rTo ^{\lambda _{X,B}}&X\pitchfork 
T(B)&\rTo ^{X\pitchfork \beta}& X\pitchfork B 
\end{diagram},\label{cotensor}\end{equation}
using the cotensorial strength $\lambda$ of $T$.

In \cite{CCGBCM},  we proved that if $T$ is a commutative monad, then 
$(A,\alpha )\pitchfork _{T}(B,\beta )$, as a subobject of $A\pitchfork 
B$, is actually a sub-{\em $T$algebra} (with the algebra structure of 
$A\pitchfork B$ given by the recipe above, with $X=A$). This in fact 
makes $\E ^{T}$ into a closed category in its own right, cf.\ Theorem 2.2 
in \cite{CCGBCM}. (It is even a {\em symmetric} closed category, in 
the sense suggested in loc.cit.\ ; this was substantiated in 
\cite{deSchipper2}, \cite{deSchipper1}.)

The notion of cotensor and $\E$-valued hom are related by an 
($\E$-strong)  
adjointness, as is well known, cf.\ \cite{Kelly} (3.42). This implies that 
for $B=(B,\beta )$ in $\E ^{T}$, we have contravariant functors
$$-\pitchfork (B,\beta ): \E \to \E^{T}$$
and
$$-\pitchfork _{T}(B,\beta ) :\E^{T}\to \E$$
which are strongly adjoint to each other on the right, so that we get 
a strong monad on $\E$, with functor part
$X\mapsto (X\pitchfork (B,\beta ))\pitchfork _{T}(B,\beta )$, or 
with slight abuse of notation
\begin{equation}\label{ddmx}X\mapsto (X\pitchfork B)\pitchfork 
_{T}B,\end{equation}
a ``restricted double dualization'' functor (terminology from 
\cite{PRSFA}). 
These double dualization monads are rarely commutative (even for 
commutative $T$); and their categories of 
algebras are often hard to analyze.

In case where $\E$ is the category of sets, and $T(X)$ is the monad 
whose algebras are boolean algebras, $(X\p 2)\p _{T}2$ is the set of 
ultrafilters on $X$, and the category of algebras for $(-\p 2)\p 
_{T}2$ is the category of  compact Hausdorff spaces (Manes' Theorem, 
cf.\ \cite{StoneSpaces} III.2.4).

If $T$ is the identity functor, and $D\in \E$ is any 
object, we have the ``plain'' double-dualization monad $(-\p D)\p D$, 
studied in detail in \cite{DD}. It is the ``full algebraic theory of 
$D$'', if we identify monads on the category of sets with infinitary 
Lawvere theories (as is done in \cite{ManesBook}, or \cite{Wraith}).

It is easy to see that $B$ itself is an algebra for the 
(unrestricted) double dualization monad $(-\p B)\p B$; the structure 
is the map $ev_{id}:(B\p B)\p B \to B$ which is ``evaluation at the identity 
map $id_{B}\in B\p B$''. In particular
\begin{equation}\label{expecx} ev_{id}\circ \delta _{B}= id_{B}.
\end{equation}

\medskip

Another significant example: if $R$ is a commutative ring object in 
$\E$, there is (under suitable completeness conditions on $\E$) a 
monad $T$ whose category of algebras are the $R$-module objects. So 
$R$ itself is a $T$-algebra (in fact, $R=T(1)$), and we have the restricted double 
dualization monad
$(-\p R)\p _{T}R$. In some examples  $(X\p R)\p _{T}R$  can be analyzed as an 
internal version of {\em distributions with compact support} on $X$ 
(distributions in the sense of Schwartz); see \cite{MR} Prop.\ II.3.6 for 
some toposes of $C^{\infty}$ spaces. An algebraic analysis is given in 
\cite{PRSFA}  for the case where $R$ is the generic 
commutative ring. 

We return to the general case of a restricted double dualization 
monad $X\mapsto (X\pitchfork B)\pitchfork 
_{T}B$, where $T$ is a strong monad on $\E$ and $B=(B,\beta )$ a 
$T$-algebra. The unit for this monad 
is denoted $\delta$, so 
$$\delta _{X}: X\to (X\pitchfork B)\pitchfork 
_{T}B.$$
Post-composing with the inclusion $(X\p B)\p _{T}B \subseteq (X\p 
B)\p B$ gives the combinator $\delta  _{X}$ considered in Section \ref{BCX}.
If $\E$ is the category of sets, it is the map which takes $x\in X$ 
to the $T$-algebra map $\delta _{X}(x):X\p B \to B$, ``evaluating at 
$x\in X$''. This ``evaluation at $x$''
is a $T$-homomorphism, thus an element in $ (X\pitchfork B)\pitchfork 
_{T}B$. In distribution terms, it is the Dirac distribution on $X$ at 
$x$, whence the notation $\delta$. The $\overline{\mu}$ of the monad $(-\p B)\p 
_{T}B$ can also ultimately be described in terms of $\delta$. We 
describe it when, and  to the extent we need it, in the proof of  
Theorem \ref{first} below.

These double-dualization monads depend on the choice of the object 
($T$-algebra) $B$. The most important case for us is where $B$ is 
$(T(1), \mu _{1})$ (later on, we shall denote this particular 
$T$-algebra by the letter $R$; it plays the role of a number line). Recall that for any $X\in \E$, the algebra $(T(X), 
\mu _{X})$ is the {\em free} $T$-algebra on $X$. In particular, $T(1)$ 
is the free algebra in one generator. 

The construction of the restricted double dualization monad 
(\ref{ddmx}) does not depend on commutativity of the given monad $T$, 
however, the following does.
So let $T$ be a {\em commutative} monad on $\E$, and let $B=(B,\beta 
)$ be a $T$-algebra. 
Then by \cite{CCGBCM}, $A\p _{T}B$ carries structure of a $T$-algebra 
whenever $A=(A,\alpha )$, $B= (B,\beta )$ are $T$-algebras. 
Therefore, the map $\delta : X \to (X\p B)\p _{T}B$ extends uniquely 
to a $T$-homomorphism on the free $T$-algebra $T(X)$, so that we have 
a canonical $T$-homomorphism
$$\tau _{X}: T(X) \rTo (X\p B)\p _{T}B.$$
Its relationship to $\int _{X}$ is made explicit in Proposition 
\ref{expadjx} 
below.
\begin{thm}\label{first}Let $T$ be a commutative monad. Then the maps 
$\tau _{X}$ form a strongly natural morphism 
$\tau :T 
\Rightarrow (-\p B)\p _{T}B$; it is a morphism of monads, and  it is 
compatible with the canonical monoidal structures on the functors in 
question.
\end{thm}
(This holds, whether one takes $\psi$ or 
$\tilde{\psi}$ as the monoidal structure  on the double dualization 
monad; and for $T$, there is 
anyway only one canonical monoidal  struture, since $T$ is assumed commutative.)

\medskip
\noindent{\bf Proof.} By construction, $\tau \circ \delta = \eta$, so 
$\tau$ is compatible with the units of the two monads in question. 
Let us prove compatibility with the $\mu$ , $\overline{\mu}$ (the 
latter being the multiplication of the double dualization monad in 
question). The 
unit and counit of the adjunction that gave rise to the monad $(-\p 
B)\p _{T}B$ are $\delta _{X}:X\to (X\p B)\p _{T}B$ in $\E$ (already 
considered), and $\epsilon _{A}:A\to (A\p _{T}B)\p B$ in $\E ^{T}$, 
given by essentially the same recipe which gave $\delta$ 
(the counit goes the ``wrong way'' because of the contravariant 
nature of the two adjoint functors in question). Since the 
multiplication of a monad arising from an adjoint pair is an instance of the counit of it, we conclude 
that $\overline{\mu}$ does indeed live in $\E ^{T}$.
Therefore, the two maps to be compared to prove that $\tau$ is 
compatible with $\mu, \overline{\mu}$ are two maps $T^{2}(X) \to (X\p 
B)\p _{T}B$ both of which are $T$-linear. By Proposition \ref{twop}, 
it therefore it suffices to 
see that they agree when precomposed with $\eta _{TX}$. Here is the 
relevant diagram:
$$\begin{diagram}TX&\rTo ^{\eta _{TX}}&T(TX) &\rTo ^{T(\tau )}&T((X\p B)\p _{T}B)& \rTo 
^{\tau} &((((X\p B)\p _{T}B)\p B)\p _{T} B)\\
&&\dTo^{\mu _{X}}&&&&\dTo_{\overline{\mu}_{X}}\\
&&TX&&\rTo _{\tau }& &(X\p B)\p _{T}B
\end{diagram}.$$
The counterclockwise composite gives $\tau _{X}$, by the unit law for 
the monad $T$. The top composite may be rewritten, using 
naturality of $\eta$, into the composite
$$\begin{diagram}TX&\rTo ^{\tau _{X}}&(X\p B)\p _{T}B&\rTo^{\eta}& 
T((X\p B)\p _{T}B)&\rTo ^{\tau }& ((((X\p B)\p 
_{T}B)\p B)\p _{T}B);
\end{diagram}$$
but  $\tau \circ \eta$ is $\delta$, by construction of $\tau$, and 
$\delta$ composed with $\overline{\mu}$ is an identity map (one of 
the monad laws for the double dualization monad here). So the 
clockwise composite likewise ends up as $\tau _{X}$.
So $\tau$ is indeed a morphism of monads. Since everything is 
compatible with the strengths, we conclude from Proposition 
\ref{montrans}  that $\tau$ also 
preserves the monoidal structure.

\medskip
\noindent{\bf Remark.} This theorem is analogous to Theorem 3.2 in 
\cite{DD}; there, however, one considers the {\em full} double dualization 
monad $(-\p B)\p B$ for an {\em unstructured} object $B$.

\medskip
The transformation $\tau$ in the Theorem is in fact an exponential 
adjoint version of the ``integral'' studied in Section \ref{BCX}:
\begin{prop}\label{expadjx}The map $\tau _{X}:T(X)\to (X\p B) \p  B$
has for its exponential adjoint the map $\int _{X,B}: T(X)\times (X\p 
B)\to B$.\end{prop}
{\bf Proof.} Since $\tau$ is the $T$-linear extension of $\delta :X\to 
(X\p B)\p B$, $\tau$ may be described explicitly as the composite
\begin{equation}\begin{diagram}TX&\rTo^{T(\delta )}&T((X\p B)\p B)&\rTo 
^{\lambda}&(X\p B)\p TB&\rTo ^{id \p \beta }&(X\p B)\p B.
\end{diagram}\label{explicitx}\end{equation}
Thus, the exponential adjoint of $\tau$ appears as the clockwise 
composite in (\ref{intxx}) composed with $\beta$.
On the other hand, if one follows the counterclockwise composite in 
(\ref{intxx}) by $\beta$, we  see from (\ref{intx}) that we have the map $\int 
_{X,B}$- 

\medskip

 Let $(B,\beta )$ be a $T$-algebra. By the Theorem, we have for each $X\in \E$ a 
map $\tau _{X}: T(X)\to (X\p B)\p _{T}B$, defined in terms of 
$\beta$, and with good properties, in 
particular, it is $T$-linear. For $X=B$, we have in particular $\tau 
_{B}:T(B)\to (B\p B)\p _{T}B$. We have also a map $ev_{id}:(B\p 
B)\p _{T}B \to B$ ``evaluation at 
$id_{B}$'', and thus get by composition a map $T(B)\to B$.
\begin{prop}\label{expectationx}The composite 
$$\begin{diagram}T(B) & \rTo ^{\tau _{B}}&(B\p B)\p 
_{T}B&\rTo^{ev_{id}}&B
\end{diagram}$$
equals $\beta : T(B)\to B$.
\end{prop}
{\bf Proof.} Both maps to be compared are $T$-linear, so it suffices 
to see that they agree when precomposed with $\eta _{B}$. We have 
$\beta \circ \eta _{B} = id_{B}$, by the unit law for $T$-algebras. 
On the other hand, $$ev_{id}\circ \tau_{B}\circ \eta_{B} = 
ev_{id}\circ \delta _{B}$$
by construction of $\tau$, and $ev_{id}\circ \delta _{B}=id_{B}$, as 
we observed above (\ref{expecx}) for unrestricted double dualization into $B$; it 
holds then, by restriction, also for the restricted double 
dualization monad.

\medskip
Theorem \ref{first} allows us to describe $\tau _{X\times Y}(\psi 
(P,Q))$ in terms of $\tau _{X}(P)$ and $\tau _{Y}(Q)$, and similarly 
for $\tilde{\psi }$; note the formal similarity with Fubini's Theorem. 
\begin{thm} Let $P\in T(X)$ 
and $Q\in T(Y)$, and let $\phi \in (X\times Y)\p B$. Then $\tau _{X\times Y}(\psi 
(P,Q))(\phi )$ and $\tau _{X\times Y}(\tilde{\psi} 
(P,Q))(\phi )$ appear as the left and right hand side, respectively, of the following 
equation (which therefore {\em holds} for all $P,Q, \phi$ when $T$ is 
assumed to be a 
commutative monad)
 \begin{equation}\int _{X}\biggl( \int _{Y}\phi (x,y) \; 
dQ(y)\biggr)\; dP(x) =\int 
_{Y}\biggl(\int _{X}\phi (x,y)\; dP(x)\biggr)\; dQ(y).
\label{fubini}\end{equation} 
\end{thm}
{\bf Proof.} We first argue that $\tau _{X\times 
Y}(\psi_{X,Y}(P,Q))(\phi ) \in B$ is given by the expression on the 
left hand side. We denote the combinators for the strong monad $S=(-\p 
B)\p_{T}B$ by $\overline{t}'$, $\overline{\psi}$, etc. Then by Theorem 
\ref{first}
$$\tau _{X\imes Y}(\psi (P,Q))= \overline{\psi} _{X,Y}(\tau 
_{X}(P),\tau _{Y}(Q)).$$
Therefore, it is a matter of analyzing $\overline{\psi}_{X,Y}$ for 
the monad $S$, and this is pure $\lambda$-calculus; in fact, $S$ can 
easily be proved to be a 
submonad of the full double dualization monad $D=(-\p B)\p B$. We claim 
that the monoidal structure $\overline{\psi}$ for this monad is 
given, in elementwise terms, as follows, for $\overline{P}\in (X\p 
B)\p B$, $\overline{Q}\in (Y\p B)\p B)$ and $\phi \in (X\imes Y)\p B$:
$$\overline{\psi}(\overline{P},\overline{Q})(\phi ) 
=\overline{P}[x\mapsto \overline{Q}[y\mapsto \phi (x,y)]].$$
This is an elementwise reformulation of the following (writing $\psi$ 
for $\overline{\psi}$ etc.:
\begin{prop} The monoidal 
structure $\psi _{X,Y}:D(X)\times D(Y) \to D(X\times Y)$ on the 
functor $D=(-\p B)\p B$ may be described 
as follows: for $P\in D(X)$ and $Q\in D(Y)$, the value of $\psi 
_{X,Y}(P,Q)$ on $\phi \in (X\times Y)\p B$
is given as the composite
\begin{equation}\label{newpsi}\begin{diagram}(X\times Y)\p B &
\rTo ^{\cong}&X\p (Y\p B)&\rTo ^{X\p Q}&X\p B&\rTo ^{P}&B.
\end{diagram}\end{equation} Similarly, the value of $\tilde{\psi}(P,Q)$ on 
$\phi$ is given as the composite
$$\begin{diagram}(X\times Y)\p B &\rTo ^{\cong}&Y\p (X\p B)&\rTo 
^{X\p P}&Y\p B&\rTo ^{Q}&B
\end{diagram}.$$
 \end{prop}
 \noindent {\bf Proof.} One may prove this by brute force, by 
$\lambda$-calculus, but since $D(D(X\times Y))$ is involved, this 
means  that a four times dualization into $B$ is involved, and this is not easy to handle; 
some ML type program on a computer would be useful here! However, we can use the 
fact that $\psi$ is characterized by being linear in the first 
variable, and to restrict along the unit (here: $\delta$) to $t''$, 
cf.\ Proposition \ref{charpsi}. So we shall prove that (\ref{newpsi}), as a function of 
$P,Q$, satisfies these two criteria. We shall be content with arguing 
elementwise (synthetically). So 
consider $P\in D(X)$ (so $P:X\p B \to B$) and $Q\in D(Y)$ (so $Q:Y\p 
B\to B$). Then (\ref{newpsi}) returns with $P$, $Q$ as input the 
composite described. We must argue that it, for fixed $Q$, depends 
linearly on $P$; recall that ``linear'' presently means 
``$D$-linear'', i.e.\ ``homomorphisms of $D$-algebras''. The
function of $P$ given by (\ref{newpsi}) is the map
$$\begin{diagram}(X\p B)\p B & \rTo ^{s\p B}& ((X\times Y)\p B)\p B
\end{diagram},$$
where $s$ is the map 
$$\begin{diagram}(X\times Y)\p B &\rTo ^{\cong}& X\p(Y\p B)&\rTo 
^{X\p Q}&X\p B.
\end{diagram}$$
Now any object of form $U\p B$ is canonically a $D$-algebra, and any 
morphism $V\p B \to U\p B$ of the form $s\p B$ (for $s:U\to V$) is a 
$D$-algebra homomorphism, since $-\p B: \E \to (\E ^{D})^{op}$ is the 
left adjoint of the two adjoint functors that together produced the 
monad $D$.

To prove the other condition, ``precomposing with $\eta$'', consider 
what happens if one puts $P = \delta_{X}(x) =$ evaluation at $x$, 
where $x\in X$ (recall that $\eta$ now 
is Dirac delta formation). Then $P$ gets replaced by $ev_{x}$, so the 
value of (\ref{newpsi}) is 
$$\begin{diagram}(X\times Y)\p B &\rTo ^{\cong}&X\p (Y\p B)
&\rTo ^{X\p Q}&X\p B&\rTo ^{ev_{x}}&B.
\end{diagram}$$ But $ev_{x}\circ (X\p Q) = Q\circ ev_{x}$, and 
precomposing $ev_{x}$ with the isomorphism $(X\times Y)\p B \cong 
X\p (Y\p B)$ yields $Q\circ (\tilde{u}_{x}\p B)$, and thus we arrive at (the value 
of $t''_{X,Y}$ at $(x,Q)$), as given  at the beginning of Section 
\ref{BCX} (replacing $T$ by $(-\p B)\p B$).

\section{Monads and actions by monoids}\label{MAMX}

Exploiting the fact that 
the functor $T$ carries two monoidal structures, we get in particular 
that $T(1)$ carries two natural 
monoid structures, namely, first,
\begin{equation}\label{m1}\begin{diagram}
T(1)\times T(1)& \rTo ^{\psi _{1,1}}&T(1\times 1) \cong T(1)
\end{diagram}\end{equation}
and, secondly, the one obtained by replacing $\psi$ with 
$\tilde{\psi}$.  They of course agree when $T$ is commutative.
The monoid multiplication $m$ (\ref{m1}) may be described equivalently as the 
composite
\begin{equation}\label{m2}\begin{diagram}
T(1)\times T(1)&\rTo ^{t' _{1,T(1)}}&T(1\times T(1))&\cong T^{2} 
(1)&\rTo ^{\mu _{1}}&T(1).
\end{diagram}
\end{equation}
This follows by recalling the construction of $\psi$ in terms of $t'$, 
$t''$ and $\mu$, and noting that $t''_{1,1}$ may be eliminated, since 
it equals the composite of 
two ``trivial'' isomorphisms $1\times T(1) \cong T(1) \cong T(1\times 1)$, cf.\ 
\cite{MSMCC} Lemma 1.8  (in fact, in the cartesian closed case, one has more generally commutativity of
$$\begin{diagram}X\times TY&\rTo ^{t''_{X,Y}}&T(X\times Y)\\
\dTo^{pr}&&\dTo _{T(pr)}\\
TY&\rTo _{=}&TY
\end{diagram}$$
where $pr$ denotes the projection). 
From either description (\ref{m1}) or (\ref{m2}) follows that the 
multiplication of $T(1)$ is 1-linear. (It is not necessarily 2-linear, 
even when $T(1)$ happens to be commutative. However, if $T$ is 
commutative, the multiplication is bilinear.)

The unit $e$ of the monoid $T(1)$ is $\eta _{1}:1\to T(1)$, also 
sometimes denoted $1$.

Any object of the form $T(X)$ carries a left action 
by $T(1)$, and also a right action by $T(1)$, the latter (which will 
be our main concern) given by
$$\begin{diagram}
T(X)\times T(1)& \rTo ^{\psi _{X,1}}&T(X\times 1) \cong T(X).
\end{diagram}$$
This action is unitary and associative, using the monoid structure on $T(1)$ given by $\psi$; 
if one prefers, one can replace simultaneously 
 $\psi _{X,1}$ and $\psi _{1,1}$ by the corresponding $\tilde{\psi}$s. (For the
 left action by $T(1)$,  one uses either $\psi$ for 
both the action and the monoid structure, or uses $\tilde{\psi} $ for 
both the action and the monoid structure.) We stick to right 
action, defined by $\psi$, as in the  displayed formula.  It is immediate 
to see that if $f:X\to Y$ is any map, then the map $T(f):T(X)\to 
T(Y)$ is equivariant for the action.

\medskip

The action of the monoid $T(1)$ on $ T(X)$  may be discussed (for 
some of its aspects) in more 
generality as follows: Let $T=(T,\eta ,\mu , t')$ be a strong monad on $\E$, and let 
$R=(R,e,m)$ be a monoid in $\E$ (with $e:1\to R$ the unit and 
$m:R\times R\to R$ the multiplication). There is an evident notion of a 
$T$-linear right action of $R$ on $T$, namely a family of unitary and 
associative actions 
(natural in $X\in \E$)
$\vdash _{X}: T(X)\times R \to T(X)$, with $\vdash _{X}$ 1-linear.

A  1-linear action  by a monoid $R$ on the monad $T$ is 
by Proposition \ref{twop} determined by its restriction (for each 
$X$) along 
$\eta _{X}\times R$, i.e.\ by  maps $\rho _{X}: X\times R \to T(X)$, 
  natural in $X$. So the unit and associativity constraints for 
the action can be encoded in terms of $\rho$. We have
\begin{prop}The 1-linear extension  of a map $\rho :X\times R \to 
T(X)$ satsfies the unit constraint iff 
$$\begin{diagram}\biggl( X\cong X\times 1& \rTo^{X\times e}&X\times R&\rTo 
^{\rho}&T(X)\biggr) = \eta _{X}
\end{diagram},$$ 
and it satisfies the 
associativity constraint iff the following diagram commutes:
\begin{equation}\label{rhox}\begin{diagram}
X\times R \times R&\rTo ^{\rho \times R}&T(X)\times R&\rTo 
^{t'_{X,R}}&T(X\times R)&\rTo ^{T(\rho )}&T^{2}(X)\\
\dTo^{X\times m}&&&&&&\dTo _{\mu _{X}}\\
X\times R&&&\rTo _{\rho}&&&T(X).
\end{diagram}\end{equation}
\end{prop}
{\bf Proof.} We leave the proof of the first assertion to the reader. 
Assume now that (\ref{rhox}) commutes. To prove that 
the action  is associative means proving equality of two maps 
$T(X)\times R \times R \to T(X)$, both of which are 1-linear. So it 
suffices to prove that these two maps agree when precomposed with 
$\eta _{X}\times R \times R$. The resulting diagram is then seen to 
be (\ref {rhox}); note that the three last arrows in the clockwise 
composite in (\ref{rhox}) is just the action, by the explicit formula 
(\ref{extx}) 
for how a map $X\times R \to T(X)$  extends to a 1-linear map 
$T(X)\times R \to T(X)$. -- We leave to the reader the proof that 
associativity of the action  implies commutativity of 
(\ref{rhox}) (and  we shall not need this implication).
 
\medskip

We return to the special case of the right action by the monoid $T(1)$ on 
$T(X)$. We denote this action simply by a dot, $P\cdot \lambda$, for $P\in T(X)$ and $\lambda \in T(1)$. We think 
of $T(1)$ as ``scalars''. 

We ask the question whether not only {\em free} $T$-algebras, but 
general $T$-algebras carry an action by $T(1)$. For this, we need 
commutativity of $T$; we have the following (which is not used in the 
sequel).

\begin{prop}\label{three}Let $T$ be a commutative monad, and let $(A,\alpha)$ be a 
$T$-algebra. Then $A$ carries a unique action by the monoid $T(1)$, 
in such a way that $\alpha :T(A)\to A$ is equivariant. The action is 
unitary and associative, and any homomorphism of $T$-algebras is 
equivariant.
\end{prop}
{\bf Proof.} By general 
monad theory, we have that $\alpha :T(A)\to A$ sits in a canonical 
absolute coequalizer diagram in $\E $
$$T^{2}(A) \tto T(A) \to A,$$
where the two parallel maps are $T(\alpha )$ and $\mu _{A}$ 
respectively. The map $T(\alpha )$ is equivariant for the action, 
without any assumptions on $T$.  We shall use commutativity of $T$ 
to prove equivariance of $\mu _{A}$. When this is established, it is 
clear that the action descends along $\alpha$  from $T(A)$ to $A$, and the rest is 
then easy. Equivariance of $\mu _{A}$ means that the right hand 
region in the following diagram commutes:
$$\begin{diagram}[nohug]TA\times 1&\rTo ^{\eta_{TA}\times \eta 
_{1}}&T^{2}A\times T1 & \rTo ^{\mu _{A}\times T1}& TA\times T1\\
\dTo^{pr}_{\cong}
&\rdTo_{\eta _{TA\times 1}}&\dTo_{\psi _{TA,1}}&&\dTo _{\psi _{A,1}}\\
TA&&T(TA \times 1)&&T(A\times 1)\\
&\rdTo_{\eta _{TA}}&\dTo _{T(pr)}^{\cong}&&\dTo _{T(pr)}^{\cong}\\
&&T^{2}A&\rTo_{\mu _{A}}&TA
\end{diagram}$$
Here $pr$ denotes the relevant projections. The triangle commutes 
because $\eta$ is a monoidal transformation. The (slanted) square 
commutes by naturality. So the counterclockwise composite equals 
$pr:TA\times 1 \to TA$. The top line composite is just $TA\times \eta 
_{1}$. The clockwise composite of the total diagram is $T(pr)\circ 
\psi _{A,1}\circ (TA\times \eta _{1})$; this, however, is again just 
$pr:TA\times 1\to TA$, by a general law for the relationship between 
$\eta$, $\psi$ and the unit isomorphisms (here the $pr$), cf.\ \cite{MSMCC}, diagram (2.3).
So we conclude that the total diagram commutes. Now, the two 
composites in the right hand region are both $T$-bilinear, because the 
$\psi$s 
are $T$-bilinear by commutativity of $T$. So to prove commutativity of 
the right hand region, it suffices to prove that it commutes after 
precomposotion with $\eta  \times \eta$, which is what the 
commutativity of the total diagram expresses. This proves the 
Proposition.

\medskip

It is easy to see that if the monad $T$ is $M\times -$ for a 
non-commutative monoid $M$ in the category of sets, then $\mu _{A}$ will 
not be equivariant; so for Proposition \ref{three}, one cannot dispense with the commutativity 
assumption for the monad $T$.

\medskip

Even though the projection $pr: A\times 1 \stackrel{\cong}{\to} A$ appears in the above 
construction and argument,  all the 
constructions and arguments work in general symmetric monoidal closed 
categories, using the unit object $I$ instead of $1$, and using the 
unit isomorphisms (part of the data of a monoidal closed category) 
$A\otimes I \cong A$ instead of $pr$. The  construction in the 
following Sections, 
however, depend in a crucial way of the assumption that our category 
is {\em cartesian} closed.

\section{Action  of functions on  distributions}\label{ABTX}
We consider a  strong (not necessarily commutative) monad $T$ on $\E$. But from now, we assume 
not only that $\E$ is  symmetric monoidal 
closed, but that it is {\em cartesian} closed (as the notation 
in the previous sections  anticipated). 
Then the object $1$ is terminal, and we have the notion of the total; 
for $P\in T(X)$, its {\em total} is $T(!)(P)\in T(1)$, where $!:X\to 
1$ is the unique such map.

\medskip 

Recall that $T(1)$ carries a canonical monoid structure, $m$,$e$, 
with $m$ defined in terms of $\psi _{1,1}$, cf.\ (\ref{m1}). 

\begin{prop}\label{prodx}Let $P\in T(X)$ and $Q\in T(Y)$. Then the total of 
$\psi_{X,Y}(P,Q)$ is the product by $m$ of the totals of $P$ 
and $Q$.
\end{prop}
This is an immediate consequence of the definition of $m$ together 
with naturality of $\psi$ with respect to the maps $!:X\to 1$ and 
$!:Y\to 1$.

\medskip
 
The space $X\p T(1)$ inherits a monoid structure $m_{X}$, $e_{X}$ in a 
standard (``pointwise'') way. 
We shall equip any free $T$-algebra $T(X)$ with a 
(right) $T$-linear action $\vdash$ by the monoid $X\p T(1)$. The construction does not depend on 
commutativity of the monad $T$.  We shall construct a map
\begin{equation}\label{actx}\begin{diagram}T(X)\times (X\p T(1)) & \rTo ^{\vdash} 
T(X)\end{diagram}\end{equation} 
It is defined as the unique 1-linear extension over $\eta _{X}\times 
T(1)$ of the following composite map:
\begin{equation}\label{primi}\begin{diagram}\rho :=&&X\times (X\p T(1))&\rTo ^{\langle pr , ev 
\rangle}&X\times T(1)& \rTo ^{t''_{X,1}}& T(X\times 1)\cong T(X).
\end{diagram}\end{equation}
Here, $pr$ denotes the projection $X\times (X\p T(1)) \to X$ to the 
first factor, and 
$ev$ denotes the evaluation map $X\times (X\p T(1)) \to T(1)$. The 
composite map displayed is actually 2-linear, so if the monad $T$ happens 
to be commutative, the 1-linear extension of it to $\vdash $ will 
 be bilinear, by Proposition \ref{fivex}. 

By Proposition \ref{precomp}, it is clear that an alternative 
description of $\rho$ is:
$$\begin{diagram}X\imes (X\p T1)&\rTo ^{\langle pr, ev \rangle}&X\times 
T1&\rTo^{\eta_{X}\times id}&TX\times T1&\rTo^{\psi_{X,1}}T(X\times 
1)\cong TX.
\end{diagram}$$

The action $\dashv$  of  $X\p T(1)$ on $T(X)$ presented here (``action 
by {\em functions} on distributions'') restricts to the 
action  of $T(1)$ on $T(X)$ (``by {\em scalars} on distributions'') 
considered in Section \ref{MAMX}, via 
the monoid map $!\p T(1): 1\p T(1) \to X\p T(1)$ induced by $! :X\to 
1$; expressed synthetically, if $\phi :X\to T(1)$ has constant value 
$\lambda \in T(1)$, then $P\vdash \phi = P\cdot \lambda$, where 
$\vdash$ denotes the action of $X\p T(1)$, and the dot denotes the 
action of $T(1)$ on 
$T(X)$.
 
\begin{thm}\label{assx}The action $\vdash :T(X)\times (X\p T(1)) \to T(X)$ is 
associative and unitary.
\end{thm}
{\bf Proof.} Our proof is not quite straightforward; there ought to 
be a better one. 
To prove the associativity assertion, we should compare two map 
$T(X)\times (X\p T(1) \times (X\p T(1) \to T(X)$ which both are 
1-linear, so it suffices to prove that their precomposite with $\eta 
_{X}\times id$ are equal. This is achieved by a contemplation of the 
following 
 diagram. (For the arrow denoted ``$\langle pr,ev \rangle $'', the middle factor $T1$ does not participate in the $\langle 
pr,ev \rangle$-formation (so elementwise, the map takes $(x,\lambda 
,\phi)$ to $(x,\lambda ,\phi (x))$); also, isomorphisms $X\times 
\cong X$ are omitted from notation.) 
$$\begin{diagram}[nohug]X\times (X\p T1)\times (X\p T1)&\rTo ^{\langle pr, ev \rangle\times 
id}&X\times T1\times (X\p T1)&\rTo ^{t_{X,1,X\p T1}}&T(X\times (X\p T1))\\
\dTo ^{X\times m_{X}}&&\dTo^{``\langle pr, ev \rangle "}&&\dTo _{T\langle pr, ev \rangle}\\
X\times (X\p T1)&&X\times T1\times T1&\rTo _{t_{X,1,T1}}&T(X\times T1)\\
&\rdTo _{\langle pr, ev \rangle}&\dTo _{X\times m}&&\dTo _{T(t''_{X,1})}\\
&&X\times T1&&T^{2}X\\
&&&\rdTo_{t''_{X,1}}&\dTo_{\mu_{X}}\\
&&&&TX;
\end{diagram}$$
the left hand region commutes,  by definition of $m_{X}$ 
in terms of $m$, and the upper right hand square is essentially just 
a twisted version of the  naturality square for $t''$ w.r.to the map 
$\langle pr, ev \rangle: X\imes (X\p T(1)) \to X\times T(1)$, recalling that th 
combinator $t$ in (\ref{mtx}) came about by a twisting of 
$t''_{X\times Z,Y}$ (here 
with $Y=1$ and $Z=T(1)$). 
The lower right hand region deserves a more detailed argument. 
Let us prove its commutativity, without using identifications 
like $X\times 1 \cong X$. Consider namely
\begin{equation}\label{starx}\begin{diagram}[nohug]X\times T1\times T1&&\rTo ^{t_{X,1,T1}}&&T(X\times 
1\times T1)\\
\dTo ^{X\times t'_{1,T1}}&&&& \\
X\times T(1\times T1)&&*&&\dTo_{T(t''_{X\times 1, 1})}\\
\dTo^{X\times T(t''_{1,1})}&&&&&\\
X\times T^{2}(1\times 1)&\rTo ^{t''_{X,T(1\times 1)}}&T(X\times 
T(1\times 1))&\rTo ^{T(t''_{X, 1\times   1})}&T^{2}(X\times 1\times 1)\\
\dTo^{X\times \mu_{1\times 1}}&&&&\dTo _{\mu _{X\times 1\times 1}}\\
X\times T(1\times 1)&&\rTo_{t''_{X,1\times 1}}&&T(X\times 
1\times 1).
\end{diagram}\end{equation}
After the identification of $1\times 1$ with $1$, the left hand 
column is ($X$ times) the defining construction of $m$; and after the 
identification of $X\times 1\times 1$, the lower right hand object is 
$T(X)$. The lower region commutes because $\mu$ is a strong natural 
transformation, thus compatible with the tensorial strengths of 
$T^{2}$ and $T$; and the upper region * is an instance of the 
generalized associativity of the tensorial strengths $t'$, $t''$. 
In more detail, writing $Y$ and $Z$ for 1, to keep them apart, 
consider

$$\begin{diagram}X\times TY\times TZ&\rTo ^{X\times 
t'_{Y,TZ}}&X\times T(Y\times TZ)&\rTo^{t''_{X,Y\times TZ}}&T(X\times 
Y \times TZ)\\
&&\dTo ^{X\times T(t''_{Y,Z})}&&\dTo_{T(X\imes t''_{Y,Z})}\\
&&X\imes T^{2}(Y\times Z)&\rTo_{t''_{X,T(Y\imes Z)}}&T(X\imes 
T(Y\times Z))\\
&&&&\dTo _{T(t''_{X,Y\times Z})}\\
&&&&T^{2}(X\imes Y\times Z).
\end{diagram}$$
The top composite is $t_{X,Y,Z}$, by (\ref{tx}). The right hand 
vertical composite is $T(t''_{X\times Y,Z})$, by the associative law 
for $t''$. So the clockwise composite in this diagram equals the 
clockwise composite of * in
(\ref{starx}) (when we put $Y=Z=1$) ; and the counterclockwise 
similarly equals the counterclockwise 
in *. This proves the associativity.

To prove the unitary law, we must prove that $id\times  e_{X}:T(X)\times 
1 \to T(X)\times (X\p T1)$ followed by 
$\vdash$ is the identity map of $T(X)$ (modulo the identification 
$T(X)\times 1 \cong T(X)$). The two maps $T(X)\times 1 \to T(X)$ to 
be compared are 1-linear, so it suffices to prove that they agree when 
precomposed with $\eta _{X}$. Consider the diagram
$$\begin{diagram}[nohug]X\times 1&\rTo ^{X\times e_{X}}&X\imes (X\p 
T1)&\rTo ^{\eta _{X}\times id}&TX\times (X\ T1)\\
&\rdTo_{X\times \eta _{1}}&\dTo_{\langle pr, ev \rangle}&& 
\dTo_{\vdash}\\
&&X\imes T1&\rTo_{t''}&T(X\imes 1)
\end{diagram}$$
The square commutes by the construction of $\vdash$, and the triangle 
commutes by the pointwise nature of $e_{X}$ in terms of $e = \eta 
_{1}$. Finally, the lower composite is $\eta _{X\times 1}$. 
After the identification of $X\imes 1$ with $X$, we thus get $\eta 
_{X}$, and this proves the unitary law.

\medskip

We next address  naturality questions for the action $\vdash$, both with respect to $X$, and with 
respect to the monad $T$. It does not immediately  make sense to ask for plain 
naturality of $\vdash$ w.r.to $X$, since the domain $T(X)\times (X\p 
T(1))$ depends both 
covariantly and contravariantly on $X$, but we do have
\begin{prop}[Frobenius reciprocity] If $f:X\to Y$ is any map, the
map $T(f):T(X)\to T(Y)$ is $Y\p T(1)$-equivariant, where  $Y\p 
T(1)$ acts on $T(X)$ via the monoid homomorphism $f^{*}:Y\p T(1) \to 
X\p  T(1)$. 
\end{prop}
Here, $f^{*}$ is short for $f\p T(1) :Y\p T(1) \to X\p T(1)$.
The statement can be expressed diagrammatically  
as commutativity of the right hand region in the diagram
$$\begin{diagram}X\times (Y\p T1)&\rTo ^{\eta _{X}\times id}&
TX \times (Y\p T1) & \rTo ^{ TX \times f^{*}}&TX 
\times (X\p T1)&\rTo ^{\vdash}&TX\\
\dTo ^{f\times id}
&&\dTo^{T(f)\times id}&&&&\dTo _{T(f)}\\
Y\times (Y\p T1)
&
\rTo _{\eta _{Y} \times id }
&TY \times (Y\p T1)&&\rTo _{\vdash }&& TY
\end{diagram}.$$
In this region, both composites are 1-linear, so as in the proof of Proposition 
\ref{three}, it 
suffices to prove commutativity of the diagram when precomposed with 
$\eta _{X}\times id$. 
Then the $\eta$s may trivially be pushed to the right, using 
naturality of $\eta$ and bifunctorality of $\times$. When the $\eta$s 
come next to the $\vdash$s, we can use the defining equations 
(\ref{primi}) 
 to eliminate $\vdash$, so that the total diagram above is rewritten 
as  
$$\begin{diagram}X\times (Y\p T1)&\rTo ^{id \times f^{*}}&X\times 
(X\p T1)&\rTo 
^{\langle pr, ev \rangle}&X\times T1&\rTo ^{t''_{X,1}}&T(X\times 
1)&\cong & T(X)\\
\dTo ^{f\times id}&&&&\dTo^{f\times id}&&\dTo _{T(f\times id)}&&\dTo 
_{T(f)}\\
Y\times (Y\p T1)&&\rTo _{\langle pr, ev \rangle}&&Y\times T1&\rTo 
_{t''_{Y,1}}&T(Y\times 1)&\cong & T(Y)
\end{diagram}.$$
Here, the left hand region commutes for pure ``$\lambda$-calculus'' 
reasons, and the rest commutes by naturality. This proves the 
Proposition.

\medskip

We next consider the naturality w.r.to morphisms of monads $\tau 
:T\to S$. This is simpler:
\begin{prop}\label{fourteenx}Let $\tau : T\Rightarrow S$ be a morphism of strong 
monads. Then $\tau _{X}: T(X)\to S(X)$ is $X\p T(1)$-equivariant, where 
$S(X)$ is equipped with action by $X\p T(1)$ via the monoid 
homomorphism $X\p  \tau _{1}:X\p T(1) \to X\p S(1)$.
\end{prop}
In diagrammatic terms, this says that the following diagram commutes:
$$\begin{diagram}TX\times (X\p T1)  &\rTo ^{\vdash}& TX\\
\dTo ^{\tau _{X}\times (X\p \tau _{1})}&&\dTo _{\tau _{X}}\\
SX\times (X\p S1)&\rTo _{\vdash}&SX
\end{diagram}.$$
For, it is standard monad theory that a monad morphism $\tau : T\Rightarrow S$ 
induces a ``forgetful'' functor $\E ^{S}\to \E ^{T}$, compatible with 
the ``underlying'' functors; and then $\tau _{X}:T(X)\to S(X)$ is a 
$T$-homomorphism; similarly for ``$T$-homomorphisms  in the first 
variable'', like the left hand vertical map in the displayed diagram. 
 Since both composites thus are $T$-linear in the first variable, it 
suffices by Proposition \ref{twop}  to see that 
we get a commutative diagram when we precompose the displayed 
diagram by $\eta _{X}\times id$, and this is straightforward.


\medskip
 We address the question of the relation between the action $\vdash 
:T(X)\times (X\p T(1))\to T(X)$ 
of functions on distributions, and $T(1)$-valued integration 
$\int_{X} : T(X)\times (X\p T(1))\to T(1)$. We
express this elementwise, and leave the diagrammatic description to 
the reader. We first note
\begin{prop}Let $P\in T(X)$ and let $\phi \in X\p T(1)$. Then
$\int _{X}\phi (x)\;dP(x)$ equals the total of $P\vdash \phi$.
\end{prop}
{\bf Proof.} We are comparing the value at $P, \phi$ of two maps 
$T(X)\times (X\p T(1)) \to T(1)$. Both are 1-linear, so it suffices to 
see that they agree when precomposed with $\eta_{X}\times id:X\times 
(X\p T(1))\to T(X)\times (X\p T(1))$. Precomposing $\int_{X}$ yields 
by (\ref{intxz}) 
the evaluation map $X\times (X\p T(1))\to T(1)$. For the other 
composite, we recall the description of $\vdash$ as 1-linear extension 
of the map $\rho: X\imes (X\p T(1))\to T(X)$ in (\ref{primi}), here 
appearing as the top composite in
$$\begin{diagram}X\times (X\p T1)&\rTo ^{\langle pr, ev \rangle}&X\times T1&\rTo 
^{t''}&T(X\times 1)&\cong&T(X)\\
\dTo^{ev}&&\dTo^{!\times id}&&\dTo ^{T(!\times id)}&&\dTo_{T(!)}\\T1&\rTo_{\cong} 
&1\times T1&\rTo _{t''}&T(1\times 1)&\cong&T(1). 
\end{diagram}$$
The clockwise composite is the total in question, the 
counterclockwise is  again the evaluation map. This proves the Proposition.

\medskip
Combining with Theorem \ref{assx}, we therefore have the following 
integration theoretic significance of the action $\vdash$; again, we 
express it in elementwise terms. The monad $T$ is assumed commutative.
\begin{thm}\label{fourx}For $P\in T(X)$ and $\phi_{1}$ and $\phi _{2}$ in $X\p 
T(1)$, we have
$$\int _{X}\phi_{1}(x) \; d(P\vdash \phi_{2})(x)=
\int _{X}(\phi_{1}\cdot \phi _{2})(x)\; dP(x).$$
\end{thm}
{\bf Proof.} The left hand side is, by the Proposition,
the total of the distribution $(P\vdash \phi _{2})\vdash 
\phi_{1}$, and the right hand side is by the Proposition the total of 
the distribution $P\vdash (\phi_{1}\cdot \phi_{2})$. The result now 
follows from the associative law (Theorem \ref{assx}) for the action of the (commutative) 
monoid $X\p T(1)$ on $T(X)$.

\medskip
This Theorem can also be obtained by using the naturality of the 
various combinators with respect to transformation of monads, as 
expressed in Proposition \ref{fourteenx}; namely, one uses the 
transformation $\tau :T \Rightarrow (-\p T(1))\p _{T}T(1)$ considered 
in Theorem \ref{first}.

\section{Tensor product and convolution}\label{TPCX}
If $P\in T(X)$ and $Q\in T(Y)$, one has $\psi (P,Q)\in T(X\times Y)$. 
This is, for classical distributions, the ``tensor product'' of the 
distributions $P$ and $Q$. One has also $\tilde{\psi}(P,Q)$, which 
agrees with $\psi (P,Q)$ if the monad is commutative. We henceforth 
stick to the commutative case.

If now $m:X\times Y \to Z$ is a map, we may form the {\em 
convolution} of $P$ and $Q$ along $m$; this is
$T(m)(\psi (P,Q)) \in T(Z)$. Thus in element-free terms, convolution formation 
along $m$ is the composite
$$\begin{diagram}T(X)\times T(Y)&\rTo ^{\psi}&T(X\times Y)&\rTo 
^{T(m)}&T(Z).\end{diagram}$$
It is $T$-bilinear.

We have encountered special cases already, namely the (right) action 
of $T(1)$ on $T(X)$, which is convolution along the isomorphism 
$X\times 1 \to X$. The multiplication making $T(1)$ into a monoid is 
the special case where $X=1$, so this multiplication is likewise a 
convolution.

The convolutions that are most important in functional analysis are 
the convolutions along the addition map $+: V\times V \to V$ for an 
abelian monoid $V$; this will be a map $*: T(V)\times T(V) \to T(V)$ 
making $T(V)$ in to an abelian semigroup. Assuming that the 
monad $T$ is of the kind studied in Section \ref{MBX} below, all objects $T(X)$ 
carry a natural addition structure $+$, and $*$ and $+$ together will 
make $T(V)$ into a commutative rig. Distributivity of $*$ over $+$ 
follows from $R$-bilinearity of $*$.

\section{Physical quantities as torsors} 

To motivate the following, consider the 1-dimensional vector space $k$ 
over a field $k$. Then a $k$-linear isomorphism $k\to k$ is 
multiplication by an invertible scalar $r\in k$, and $r$ in 
fact defines a natural isomorphism $\rho : T\Rightarrow T$, where $T$ 
is the free-vector space monad, namely $\rho _{X}:T(X)\to T(X)$ is the 
homothety ``multiplication by $r$''. This transformation is compatible 
with the $\mu$ of the monad, since each instance of $\rho$ is a 
linear map; but it is not compatible with $\eta$, since 
$\rho_{1}(1)=r\in k$ 
is not necessarily $1\in k$.

\begin{prop}\label{ccx}Let $\rho : T\Rightarrow S$ be a strong natural transformation 
between endofunctors on $\E$. Then for any pair of objects
 $X,Y$, the following diagram commutes:
$$\begin{diagram}TX\times (X\p TY)&\rTo &T^{2}Y\\
\dTo^{\rho_{X}\times (X\p \rho_{Y})}&&\dTo _{\rho^{2}_{Y}}\\
SX\times (X\p SY)&\rTo &S^{2}Y
\end{diagram}$$ where the horizontal maps  are ``strength in the 
righthand factor, followed by 
evaluation'', 
 and where natural transformation  $\rho ^{2}$ denotes the  natural 
transformation $T^{2}\Rightarrow S^{2}$ derived from $\rho$.\end{prop}
Thus the top arrow is
$$\begin{diagram}TX\imes (X\p TY)&\rTo ^{TX\times st^{T}}&TX\times (TX\p 
T^{2}Y)&\rTo ^{ev}&T^{2}Y,
\end{diagram}$$
where $st^{T}$ is the strength (enrichment) of $T$; similarly  the 
bottom  arrow one is obtained from the strength $st^{S}$ of $S$.
The natural transformation $\rho^{2}$ is more explicitly given by 
$\rho ^{2}_{Y}= 
 S(\rho _{Y})\circ \rho_{TY}$. The proof of this Proposition is in 
principle elementary;
it uses that $\rho$ is a 
{\em strong} natural transformation, which means in particular (cf.\ 
\cite{Borceux} II.6.2.8) that diagrams of the form
$$\begin{diagram}X\p Y&\rTo ^{st_{T}}&TX \p TY\\
\dTo^{st_{S}}&&\dTo _{TX\p \rho _{Y}}\\
SX\p SY &\rTo_{\rho_{X}\p Y}&TX\p SY
\end{diagram}$$
commute (the equivalence of this notion of strong natural 
transformation with the one of (\ref{strx1}) is proved in 
\cite{MSMCC} Lemma 1.1).

\medskip

Let $B=(B,\beta )$ be a $T$-algebra. Recall from Section \ref{BCX} that we have a map
$$\int _{X}:T(X)\times (X\p B )\to B.$$
 Similarly for $S$. 
Inspecting the explicit construction (\ref{intbxx}) (with $(T(1); \mu 
_{1}$) for $(B,\beta)$), we note that the construction depends on 
$\mu$, but it does not 
depend on $\eta$. Therefore, the following is not surprising:  
\begin{prop}Let $T$ and $S$ be strong monads on $\E$, and let $\rho :T\Rightarrow S$ 
be a strong natural transformation, compatible with the $\mu$s, but 
not necessarily with the $\eta$s. Then the $\int _{X}$-formation for the monads $T$ and $S$ is 
compatible with $\rho$, in the sense that the following diagram 
commutes:
$$\begin{diagram}T(X)\times (X\p T(1))&\rTo ^{\int_{X}}&T(1)\\
\dTo ^{\rho _{X}\times (X\p \rho _{1})}&&\dTo_{\rho _{1}}\\
S(X)\times (X\p S(1))&\rTo _{\int_{X}}&S(1).
\end{diagram}$$
\end{prop}
{\bf Proof.} Use the explicit form 
(\ref{intbxx}) for the $\int _{X}$ in question; then the desired 
diagram comes about from the diagram in Proposition \ref{ccx} by 
putting $Y=1$, and concatenating it with the commutative square 
expressing compatibility of $\rho$ with the $\mu$s of the monads:
$$\rho _{1}\circ \mu ^{T}_{1} = \mu ^{S}_{1}\circ \rho^{2}_{1},$$
where $\mu ^{T}$ and $\mu^{S}$ denote the multiplication of the monads $T$ and $S$, 
respectively.

\medskip

Let $T$ be a commutative monad on $\E$. Consider 
another strong endofunctor $M$ on $\E$, equipped with an action $\nu$ 
by $T$,
$$\nu : T(M(X))\to M(X)$$
natural in $X$, and with $\nu$ satisfying a unitary and associative law. Then 
every $M(X)$ is a $T$-algebra by virtue of $\nu _{X}:T(M(X))\to 
M(X)$, and morphisms of the form $M(f)$ are $T$-linear.
Let $M$ and $M'$ be strong endofunctors equipped with such 
$T$-actions.
There is an evident notion of when a strong natural transformation 
$\lambda  :M\Rightarrow M'$ is compatible with the $T$-actions, so we 
have a category of $T$-actions. The endofunctor $T$ itself is an 
object in this category, by virtue of $\mu$. We say that $M$ is a 
{\em $T$-torsor} if it is isomorphic to $T$ in the category of 
$T$-actions. Note that no particular such isomorphism is chosen; this 
is just like a 1-dimensional vector space over $k$: it is isomorphic 
to $k$, but no particular isomorphism is chosen.

Our contention is that the category of $T$-torsors is a mathematical  model of 
(not necessarily pure) quantities of type $T$ (which is the 
corresponding pure quantity). Thus if $T$ 
is the free ${\mathbb R}$-vector space monad, the functor $M$ which 
to a space $X\in \E$ associates the space of distributions of 
electric charges 
over $X$, is a $T$-torsor.

\medskip

The following Proposition expresses that isomorphisms of actions $\lambda: T\cong M$ are 
determined  by $\lambda _{1}:T(1)\to M(1)$; in the example, the 
latter data means: choosing a {\em unit} of electric charge.
\begin{prop}If $g$ and $h:T\Rightarrow M$ are isomorphisms of $T$-actions, and if 
$g_{1}=h_{1}:T(1)\to M(1)$, then $g=h$.\end{prop}
{\bf Proof.} By replacing $h$ by its inverse $M\to T$, it is clear 
that it suffices to prove that if $\rho :T\to T$ is an isomorphism of 
$T$-actions, and $\rho_{1}= id_{T(1)}$, then $\rho$ is the identity 
transformation. As a morphism of $T$-actions, $\rho$ is in particular 
a {\em strong} natural transformation, which implies that right hand  
square in the 
following diagram commutes for any $X\in \E$; the left hand square commutes 
by assumption on $\rho _{1}$:
$$\begin{diagram}X\times 1&\rTo^{X\times \eta_{1}}&X\times 
T(1)&\rTo^{t''}&T(X\times 1)\\
\dTo^{=}&&\dTo^{X\times \rho_{1}}&&\dTo_{\rho_{X\times 1}}\\
X\times 1&\rTo_{X\imes \eta_{1}}&X\imes T(1)&\rTo _{t''}&T(X\imes 1)
\end{diagram}$$ 
Now both the horizontal composites are $\eta _{X\times 1}$, by general 
theory of tensorial strengths. Also $\rho _{X\times 1}$ is 
$T$-linear. Then uniqueness of $T$-linear extensions over 
$\eta_{X\times 1}$ 
implies that the right hand vertical map is the identity map. Using the natural 
identification of $X\times 
1$ with $X$, we then also get that $\rho_{X}$ is the identity map of 
$T(X)$.

\section{Monads and biproducts}\label{MBX}

Let $T$ be a commutative monad. We  summarize some of the relations between 
the covariant functor 
$T:\E \to \E$, and the contravariant $-\p T(1): \E \to \E$.
The latter is actually  valued in the category of commutative monoids 
in $\E$.

\begin{itemize}
\item There is a $T$-bilinear pairing $T(X) \times (X\p T(1)) \to 
T(1)$, namely the exponential adjoint $\int _{X}$ of the map $\tau _{X}:T(X)\to 
(X\p T(1))\p _{T}T(1)$.

\item There is an associative and unitary $T$-bilinear action 
$\vdash$ of $X\p T(1)$ on $T(X)$; it satisfies a ``Frobenius 
reciprocity'' naturality condition.
\end{itemize}

These are two of the axioms laid down by Lawvere in his description 
of relations between extensive quantities and intensive quantities, 
except for the $T$-bilinearity, cf.\ e.g.\  \cite{ACFPE}, Lecture IV; in Lawvere's axiomatics, one deals 
rather with bilinearity in the sense of an additive structure.

We shall in the present Section describe a simple categorical 
property of the monad $T$, which will guarantee that ``$T$-linearity 
implies additivity'', even ``$R$-linearity'' in the sense of a rig $R 
\in \E$ (``rig"= commutative semiring), namely $R=T(1)$. This 
condition will in fact imply that $\E^{T}$ is a ``linear category''.

We begin with some standard general category theory,
namely a monad $T=(T,\eta , \mu)$ on a category which has finite 
products and finite coproducts. (No distributivity is assumed.)
So $\E$  has an initial 
object $\emptyset$. If $T(\emptyset )\in \E$ is a terminal object,  
then the object $(T(\emptyset ), \mu 
_{\emptyset})$ is a zero object in $\E ^{T}$, i.e.\ it is both 
initial and terminal. It is initial because $T$, as a functor $\E \to 
\E^{T}$, is a left adjoint, hence preserves initials; and since 
$T(\emptyset )=1$, it is also terminal (the terminal object in 
$\E^{T}$ being $1\in \E$, equipped with the unique map $T(1)\to 1$ as 
structure). This zero object in $\E ^{T}$ we denote $0$. Existence 
of a zero object in a category implies that the category has 
distinguished zero maps $0_{A,B}:A\to B$ between any two objects $A$ 
and $B$, namely the unique map $A\to B$ which factors through $0$. 
For $\E^{T}$, we can even talk about the zero map $0_{X,B}:X\to B$, 
where $X\in \E$ and $B=(B, \beta )\in \E^{T}$, namely $\eta _{X}$ 
followed by the zero map $0_{T(X),B}:T(X)\to B$. We have a canonical 
map
$X+Y\to T(X)\times T(Y)$: the composite $X \to  X+Y \to T(X)\times 
T(Y)$ is $(\eta _{X},0_{X,T(Y)})$ (here, the first map is the coproduct inclusion 
map ). Similarly, we have a canonical map $Y \to T(X)\times T(Y)$. 
Using the universal property of coproducts, we thus get a canonical 
map $\phi _{X,Y}: X+Y\to 
T(X)\times T(Y)$. It extends uniquely over $\eta _{X+Y}:X+Y \to 
T(X+Y)$ to a $T$-linear map 
$$\Phi _{X,Y}: T(X+Y)\to T(X)\times T(Y),$$
and $\Phi$ is natural in $X$ and in $Y$.
We say that {\em $T:\E \to \E^{T}$  takes binary coproducts to 
products} if $\Phi _{X,Y}$  is an isomorphism (in $\E$ or equivalently in 
$\E^{T}$) for all $X$, $Y$ in $\E$ . Note that the definition presupposed that $T(\emptyset 
)=1$; it is the zero object in $\E^{T}$, so that if $T$ takes binary 
coproducts to products, it in fact takes finite coproducts to 
products, in a similar sense. So we can also use the phrase ``{\em $T$ 
takes finite coproducts to products}'' for this property of $T$. 

We define an ``addition'' map in $\E^{T}$
; it is a map $+: T(X)\times T(X) $ to $T(X)$, namely the composite
$$\begin{diagram}T(X)\times T(X)&\rTo ^{\Phi _{X,Y}^{-1}}& 
T(X+X)&\rTo ^{T(\nabla )}&T(X)\end{diagram}$$
where $\nabla : X+X\to X$ is the codiagonal. 
So in particular, if $in_{i}$ denotes the $i$th inclusion ($i=1,2$) of 
$X$ into $X+X$, we have
\begin{equation}\label{obv}\begin{diagram}
id_{TX} =&TX&\rTo ^{T(in_{i})}&T(X+X)&\rTo^{\Phi _{X,X}}&TX\times TX&\rTo 
^{+}&TX.
\end{diagram}\end{equation}
Under the identification $T(X)\cong T(X+\emptyset ) \cong T(X)\times 
1$, the equation (\ref{obv}) can also be read: $ T(!):T(\emptyset)\to 
T(X)$ is right unit for $+$, and similarly one gets that it is a left 
unit. 

We leave to the reader the easy proof of associativity and 
commutativity of the map 
$+:T(X)\times T(X) \to T(X)$. It follows that $T(X)$ acquires 
structure of an abelian monoid in $\E ^{T}$ (and also in $\E$).

\begin{prop}\label{linx}Every $T$-linear map $T(X)\to T(Y)$ is compatible with 
the abelian monoid structure.
\end{prop}
{\bf Proof.} This means that we should prove commutativity of the 
square * in the following diagram
$$\begin{diagram}T(X+X)&\rTo ^{\Phi }&T(X)\times T(X)&\rTo^{f\times 
f}&T(Y)\times T(Y)\\
&&\dTo^{+}&*&\dTo_{+}\\
&&T(X)&\rTo_{f}&T(Y)
\end{diagram}$$
for $f$ any $T$-linear map; so $f$ is not necessarily of the form 
$T(g)$, but it has the property that it preserves 0. To prove 
commutativity of the diagram *, it suffices to precompose with  the 
linear isomorphism $\Phi$. Now the two maps to be compared are both 
$T$-linear, and $T(X+X)$ is a coproduct in $\E ^{T}$, so it suffices to see that their composite with the 
inclusion $T(in_{i}): T(X)\to T(X+X)$ (where $i=1$ or $=2$) are equal.
Now $(f\times f)\circ\Phi \circ T(in_{i})$ is seen to be $f$, using 
that 0 is neutral for the addition.

\medskip

Recall that we have the $T$-bilinear action $T(X)\times T(1) \to 
T(X)$. It follows from the Proposition that it is additive in 
each variable separately. 

We have in particular the $T$-bilinear commutative multiplication 
$m:T(1)\times T(1)\to T(1)$, likewise bi-additive, 
$m(x+y,z)=m(x,z)+m(y,z)$, or in the notation one also wants to use,
$$(x+y)\cdot z=x\cdot z + y\cdot z,$$
 so that $T(1)$ carries structure of a rig (= commutative 
semiring). This rig we also denote $R$. The category of modules over 
a rig $R$ is the category of abelian monoids equipped with a bi-additive 
action by $R$, and maps which preserve the addition and the action. 
We may summarize:

\begin{prop}\label{modulex}Each $T(X)$ is a module over the rig $R=T(1)$; each 
$T$-linear map $T(X) \to T(Y)$ is an $R$-module morphism.
\end{prop}
It is more generally  true that  $T$-linear maps $A\to B$ (for $A$ and 
$B \in \E ^{T}$) are $R$-module maps. We shall not use this fact; it 
is proved in analogy with the proof of Proposition \ref{three}.

\medskip

Let us finally note

\begin{prop}\begin{sloppypar}If $T$ takes finite coproducts to products, then so does 
the associated Schwartz monad $S= (-\p T(1))\p _{T}T(1)$.
\end{sloppypar}
\end{prop}
{\bf Proof} (sketch). We have 
\begin{equation*}
\begin{split}(\emptyset \p T(1))\p _{T}T(1)&\cong 1\p _{T}T(1) \cong 1
\end{split}\end{equation*}
the last isomorphism because $1 =0$ is an initial $T$-algebra. 
Similarly,
\begin{equation*}
\begin{split}
S(X)\times S(Y)&= [(X\p T(1))\p _{T}T(1)]\times [(Y\p T(1))\p _{T}T(1)]\\
&= [(X\p T(1))\oplus (Y\p T(1))]\p _{T} T(1)\\
\intertext{because $\oplus$ is coproduct in $\E ^{T}$,}
&=[(X+Y)\p T(1)]\p _{T}T(1)\\
\intertext{because $\oplus$ is product in $\E$}
&=S(X+Y).
\end{split}\end{equation*}

\section{Expectation and other moments}\label{EOMX}
We consider now a  commutative monad $T=(T,\eta ,\mu , t')$ on 
$\E$ (a CCC with  coproducts and finite inverse limits), such 
that $(T,\eta ,\mu)$ takes finite coproducts to products. Thus $\E 
^{T}$ is a semi-additive category with biproducts $\oplus$; all 
its objects are modules over the rig $R=T(1)$, and all morphisms are 
$R$-linear (as well as $T$-linear, of course). We call $R$ the rig of 
{\em scalars}.

Talking synthetically, we  call the elements of $T(X)$ {\em 
concrete} distributions on $X$. We also have the object $S(X)=(X\p R)\p 
_{T}R$ of Schwartz distributions on $X$, i.e.\ $T$-linear functionals $X\p R 
\to R$; and we have the map 
$\tau _{X}:T(X)\to S(X)$ taking  concrete distributions on $X$ to 
such functionals.
We have by Proposition \ref{expadjx} that 
$\int _{X}\phi (x) \; dP(x)$
is the value of the functional $\tau _{X}(P): X\p R \to R$ on $\phi 
:X \to R$ ($\phi$ a ``test function'' on $X$, in Schwartz 
terminology).
The {\em total} of $P$ is the $T(X\to 1)(P) \in T(1)=R$, and may be 
written as $\int _{X}1_{X}\; dP(x)$ where $1_{X}:X \to R$ is 
the constant function with  the multiplicative unit $1 = e = \eta _{1}$ of $R$. 

For  concrete distributions $P$ on 
the space $R$ itself, (so $P\in T(R)$) there are other  
characteristic scalars, called ``moments'', namely some of the values 
of the functional $\tau _{R}(P): R\p R \to R$ on some particular 
functions $R\to R$.
In view of the universal role which the identity map has in e.g.\ 
Yoneda's Lemma, it is no surprise that the value of this functional 
on $id_{R}$ plays a particular role. It is the {\em expectation} of 
$P$, denoted $E(P)\in R$,
$$E(P) := \int _{R}x\; dP(x).$$
We have
\begin{prop}\label{ppx}Let $P\in T(R)$. Then any value of the functional $\tau _{R}(P): (R\p R) \to R$ 
is the expectation of some $P'\in T(R)$:
for any $\phi \in R\p R$ and $P\in T(R)$,
$$\tau _{R}(P)(\phi )= E(T(\phi )(P)).$$
\end{prop}
{\bf Proof.} By naturality of $\tau$ with respect to $\phi :R\to R$,
$$\tau _{R}\circ T(\phi ) = [(\phi \p R)\p _{T}R]\circ \tau _{R} .$$ When 
postcomposed with $ev _{id}: (R\p R)\p _{T}T \to R$, the left hand 
side gives $E(T(\phi )(P))$, the right hand side gives $\tau _{R}(P) 
(\phi )$, because $ev _{id}\circ ((\phi \p R)\p _{T}R)= ev _{\phi}$.

\medskip

Note  that for any $T$-algebra $B=(B, \beta )$, and $P  
\in T(B)$, we have $E(P)= \int _{B}x \; dP(x) = \beta (P)$; this is just a 
reformulation of (\ref{intelxx}).

\medskip

Since $R$ is a rig, we have for each natural number $n$ a map $R\to R$, 
elementwise described  by $x\mapsto x^{n}$. The {\em $n$th moment} 
$\alpha _{n}(P)$ of 
$P\in T(R)$ is defined as 
$\int _{R}x^{n}\; dP(x)$, thus $\alpha _{0}(P)$ is the total of $P$, and 
$\alpha _{1}(P)$ the expectation of $P$. Note that $\alpha 
_{1}(P)=E(P) 
=\mu_{1}(P)$, where $\mu _{1}:T^{2}(1)\to T(1)=R$ comes from the 
monad-multiplication $\mu :T^{2}\Rightarrow T(1)=R$.

\medskip
In \cite{Cramer} (5.5.6), one finds the formula $E\{ X+Y\} = E\{ X\} +E\{ 
Y\}$ where $X, Y$ is a joint distribution of two simultaneous random variables, valued in $R$. 
The formula looks 
deceptively just like it were a consequence of linearity of 
$E:T(R)\to R$ ($= \mu _{1}:T^{2}(1)\to T(1)$); but recall that  $X,Y$ is not a pair of 
distributions; rather, it is meant to denote a {\em simultaneous} distribution, i.e.\ an 
element $P \in T(R\times R)$, and $X+Y$ refers to the distribution $\in 
T(R)$ obtained by applying $T(+):T(R\times R) \to T(R)$ to $P$. So 
the formula is not a simple linearity. It is rather a formulation of 
the following:
\begin{prop}The following diagram * commutes:
$$\begin{diagram}T^{2}(2)&\rTo ^{T(\Phi )}&T(R\times R)&\rTo 
^{T(+)}&T(R)=T^{2}(1)\\
\dTo ^{\mu _{2}}&&\dTo^{\beta}&*&\dTo \mu _{1}\\
T(2)&\rTo _{\Phi}&R\times R&\rTo _{+}&R=T(1)
\end{diagram}$$
where $\beta$ is the coordinatewise $T$-algebra structure on $R\times 
R$.
\end{prop}
{\bf Proof.} Write $T(1)$ for $R$, and write $1+1$ for $2$, and let 
$\Phi$ be the comparison isomorphism, expressing that $T$ takes 
finite coproducts to products. Then the left hand square commutes, 
since $\Phi$ is $T$-linear, 
and the outer diagram commutes by naturality of $\mu$ with respect to 
the map $\nabla :2\to 1$. (Here, of course, $\nabla$ is the unique 
map $2\to 1$, but we write it for systematic reasons; in fact, the 
Proposition and the proof immediately generalizes  when $R$ is 
replaced by $R^{n}$, in which case $\nabla :2n\to n$ is not so 
trivial.)

\medskip

For comparison with the quoted formula from \cite{Cramer}, if  $X,Y$ denotes $P$, the  clockwise composite takes $P$ to $E\{ X+Y\}$, 
and the counterclockwise takes it into $E\{X \} + E\{Y \}$.

\medskip

If $P\in T(R)$ has total 1, the physical significance of $E(P)\in R$ 
is ``center of gravity'' of $P$ (thinking of $P$ as a mass 
distribution). However, physically it is clear that the center of 
gravity of a mass distribution on the line $R$ does  not the depend 
on the location 
of the origin $0\in R$, but only of the {\em affine} structure of 
$R$, in other words, it is invariant under affine maps $R\to R$. 
Here, we may take ``affine map $R\to R$'' to mean maps of the form 
$x\mapsto a\cdot x + b$ where $a$ and $b$ are scalars $\in R$.

\begin{prop}Let $P\in T(R)$ have total 1. Then for any affine $\phi 
:R\to R$,
$\phi (E(P)) = E(T(\phi )(P))$.
\end{prop}
{\bf Proof.} We may write $\phi \in R\p R$ as a linear combination of 
the identity map $id: R\to R$, and $1:R\to R$ (the map with constant 
value $1\in R$), $\phi (x)= a\cdot x + b$. By Proposition \ref{ppx},
we have
$$E(T(\phi )(P)) =\tau _{R}(P) (\phi )=  \tau _{R}(P)(a\cdot id + 
b\cdot 1).$$
Then since $\tau _{R}(P)$ is $T$-linear, 
it is $R$-linear (Proposition \ref{modulex}), so we may continue the equation
$$= a\cdot \tau_{R}(P)(id) + b\cdot \tau _{R}(P) (1),$$
which is $a\cdot E(P) + b\cdot 1$, the last term since $P$ has total 
1.

\medskip

The notion of moments  make sense not only for distributions on $R 
=T(1)$, but for instance also for distributions on $R^{2}=T(2)$.
Thus if $P\in T(2)$, we have for any $\phi:R^{2}\to R$ the
scalar $\int _{R^{2}}\phi (z)\; dP(z)$. Since the dummy variable $z$ 
here ranges over $R^{2}$, it is more natural to write it $z=(x,y)$, 
where $x$ and $y$ range over $R$, and thus the scalar in question is 
written
$\int _{R^{2}}\phi (x,y)\; dP(x,y)$. The {\em mixed second order 
momemt} of $P$ is the scalar obtained by taking $\phi$ to be the 
multiplication map $R\times R \to R$, so is $\int _{R^{2}}x\cdot y \; 
dP(x,y)$. It is in terms of this that one can define the {\em 
correlation} coefficient of $P$.

\section{Examples.}\label{exax}

The simplest example is where $\E$ is the category of sets (strength is automatic here), 
and $T$ 
is the free-commutative-monoid monad. This is related to the notion 
of ``multiset'', since $T(X)$ also may be seen as the set of 
multi-subsets of $X$; an element of $T(X)$ consists in an assignement 
$P$ 
of multiplicities $\{n(x)\in {\mathbb N}\mid x\in X \}$, with $n_{x}=0$ for all 
but a finite number of $x$s, ``$P$ is if compact support''. 
Then $T(1) ={\mathbb N}$, and $X\p T(1)$ 
is the set of assignements $\phi$ of 
multiplicities $\{n(x)\in {\mathbb N} \mid x\in X \}$, but without the 
requirement of compact support. Consider $X=T(1)={\mathbb N}$. One 
can easily see that $T({\mathbb N})$ may be identified with the set of 
polynomials in one variable with coefficients from ${\mathbb N}$, and 
then convolution along the addition map ${\mathbb N}\times {\mathbb N} 
\to {\mathbb N}$ becomes identified with multiplication of 
polynomials. (Similarly for finite products ${\mathbb N}^{k}$.)

\medskip 
An example where the conceptual machinery (strength) has to be 
brought in explicitly is the following, which was one of the 
motivations for the present research:
Consider the category $\E$ of convenient vector spaces and the smooth 
(i.e.\ $C^{\infty}$  but not 
necessarily linear) maps in between. It is a cartesian closed 
category, cf.\ \cite{KM} 
and \cite{BET}, and there exists free vector spaces (${\mathbb 
R}$-modules) in it, hence a commutative monad $T$. The category $\E$ does (probably) not have 
equalizers, at least it is clear that the zero set of a nonlinear 
map, say 
$V\to {\mathbb R}$, does not have a natural vector space structure. 
On the other hand, the equalizer of two parallel {\em linear} maps in 
$\E$ does exist. 
The following piece of general 
theory shows that therefore $\E$ has {\em enough} equalizers to form 
the subobject $A\p _{T}B \subseteq A\p B$, which was crucial in the 
construction of restricted double dualization monads (as in Section 
\ref{MDD}), and in 
\cite{CCGBCM}.

We recall from \cite{Bunge}, or \cite{CCGBCM} the two parallel maps whose 
equalizer, if it exists, gives $A\p _{T}B\subseteq A\p B$, (where 
$A=(A,\alpha )$ and $B=(B,\beta )$ are two $T$-algebras). The two maps 
$A\p B \to T(A)\p B$ are $\alpha \p B$, on the one hand, and the 
composite
\begin{equation}\label{equx}\begin{diagram}
A\p B&\rTo ^{st}&T(A)\p T(B)&\rTo^{T(A)\p \beta}&T(A)\p B.
\end{diagram}\end{equation}
The map $\alpha \p B$ is clearly $T$-linear. For the map 
(\ref{equx}), this is not immediately clear; in fact, it depends on 
commutativity of the monad $T$:
\begin{prop}Let $T$ be a commutative monad, and $A=(A,\alpha )$ and 
$B=(B,\beta )$ two $T$-algebras. Then the composite (\ref{equx}) is 
$T$-linear. 
\end{prop}
{\bf Proof.} In the diagram
$$\begin{diagram}T(A\p B)&\rTo ^{T(st)}&T(TA\p TB)&\rTo ^{T(id\p 
\beta )}&T(TA\p B)\\
\dTo^{\lambda}&&\dTo_{\lambda}&&\dTo_{\lambda}\\
A\p TB&\rTo _{st}&TA\p T^{2}B&\rTo_{id\p T\beta}&TA\p TB\\
\dTo ^{id\p \beta}&&\dTo^{id\p \mu}\dTo _{id\p T\beta}&&\dTo _{id\p \beta}\\
A\p B&\rTo_{st}&TA\p TB&\rTo_{id \p \beta}&TA\p B,
\end{diagram}$$ the vertical outer edges are the $T$-algebra 
structures on $A\p B$ and $TA\p B$, respectively, expressed in terms 
of the cotensorial strength $\lambda$. We are thus required to prove 
that the outer square commutes.Three of the inner squares commute for 
obvious reasons (ignore for the moment the arrow $id\p \mu$), 
but the upper left square does not. Now the associative law for the 
structure $\beta$ allows us to replace the ``doubled'' arrow $id\p 
T\beta $ with $id\p \mu$. But for a commutative monad $T$, the upper 
left hand square postcomposed with $id\p \mu$ commutes; this 
condition is in fact equivalent to commutativity of $T$, as stated in 
\cite{CCGBCM} Definition 2.1 (in loc.cit., it is presented as an 
alternative equivalent 
 definition of commutativity of $T$ in terms of the cotensorial strength 
$\lambda$). From this follows that the outer diagram above commutes, 
and this proves the Proposition.

\medskip

A concrete description of the monad $T$ for vector spaces over 
${\mathbb R}$ in this category was given in \cite{BET}. The authors in fact 
prove that it is``carved it out'' by topological means from a Schwartz type double 
dualization monad, described also in \cite{KM}. They provide a 
categorical study of this monad along different lines than ours, 
namely in terms of an ``exponential modality'' (essentially the 
comonad $!$ considered in linear logic).

\section{Probability theory}
To justify some measure- and probability- theoretic  terminology, one may think of 
an ele\-ment $X\p 
T(1)$ not just as a ``test function'', in the sense of 
Schwartz distribution theory, but as a generalized ``measurable 
subset'' of $X$, or as a generalized ``event'' in the ``outcome 
space'' 
$X$. The connection is that a subset $X' \subset X$ (for suitable $\E$ 
and suitable $T$) gives rise to a  function $X\to T(1)$, namely the 
characteristic function (whose value is 1 for $x\in X'$, and 0 else). 
Like for 
$X\p T(1)$, the set of subsets of $X$ depends contravariantly on $X$, 
via inverse image formation. Instead of the $T$-linearity requirement for 
Schwartz distributions $X\p T(1) \to T(1)$, there are other well 
known algebraic requirements for measures, viewed as functions from 
the boolean algebra of subsets of $X$ to the rig $T(1)$. This shall 
not concern us in detail here; the observation is just that test 
functions on $X$ may be viewed as generalized measurable 
subsets/events in $X$, and thereby it gives us access to terminology and 
notions
borrowed from measure theory or probability theory. We already 
anticipated this import of terminology when we, for $P\in T(R)$, used the word 
``expectation of $P$'' for $\int _{R}x\; dP(x)$. 

\medskip

A strong monad $T$ on a CCC $\E$ is called {\em affine} if $T(1)=1$. 
For algebraic theories, this was introduced in \cite{Wraith}. For 
strong  monads, it was proved in \cite{BCCM} that this is equivalent 
to the assertion that for all $X,Y$, the map
$\psi _{X,Y}:T(X)\times T(Y) \to T(X\imes Y)$ is split monic with 
$(T(pr_{1}),T(pr_{2})):T(X\times Y) \to T(X)\times T(Y)$ a retraction. 
In \cite{Lindner}, it was proved that if $\E$ has finite limits, any 
 commutative monad $T$ has a maximal affine submonad $T_{0}$, the 
``affine part of $T$''. It is likewise a commutative monad.
Speaking in 
elementwise terms, $T_{0}(X)$ consists of those concrete 
distributions whose total is $1\in T(1)$. We consider in the 
following a commutative monad $T$ and its affine part $T_{0}$. 

Probability distributions  have by definition total $1\in R$, (recall 
that $R$ denotes the rig $T(1)$) and  
take values in the interval from 0 to 1. We do not in the present 
article consider any order relation on $R$, so there is no ``interval 
from 0 to 1''; so we are stretching 
terminology a  bit when we use the word ``probability distribution on $X$'' for the 
elements of $T_{0}(X)$, but we shall do so. So a ``probability 
distribution'' is here just a concrete distribution $P\in T(X)$ with 
total 1, or in the notation from Section \ref{BCX}, 
$$\int _{X}1_{X}\; dP(x) =1$$
where $1_{X}:X\to R$ is the function with constant value $1\in 
R$. Since the object $1$ is terminal, it is clear that for any $f:X\to Y$, 
if $P\in T(X)$ is a probability distribution, then so is $T(f)(P)\in 
T(Y)$. (Alternatively: $T_{0}$ is a subfunctor 
of $T$.)

 If $P\in T_{0}(X)$ and $Q\in 
T_{0}(Y)$, then $\psi (P,Q)\in T_{0}(X\times Y)$, cf.\ Proposition 
\ref{prodx}; this 
also follows since the inclusion of strong monads $T_{0}\subseteq T$ is 
compatible with the monoidal structure $\psi$.
From this in turn follows that e.g.\  probability 
distributions are stable under convolution.

The assertion that $\psi _{X,Y}$ for the monad $T_{0}$ is split monic, 
quoted above, may in termino\-logy from probability theory be rendered: 
``the distribution for independent random variables may be 
reconstructed from  marginal distributions''; recall that if $Q\in 
T(X\times Y)$, then its marginal distributions are $T(pr_{i})(Q)$ 
($i=1,2$). If $Q$ is a probability distribution, then so are its 
marginal distributions.

\medskip

The subobject $T_{0}(X)\subseteq T(X)$ is clearly not stable under 
the multiplication  by scalars $\lambda \in R$; in fact, 
formation of totals is the map $T(!):T(X)\to T(1)=R$, hence is 
$T$-linear, and therefore commutes with multiplication by scalars.
In particular,  $T_{0}(X)\subseteq T(X)$ is not stable under multiplication $\vdash$ by 
functions $\phi \in X\p R$. However, this multiplication still plays a role in the 
formulation of probability theory presented here. Let $P\in 
T_{0}(X)$, and let $\phi \in X\p R$ be such that $\lambda:=\int _{X}\phi (x)\; 
dP(x)$ is invertible in the multiplicative monoid of $R$. Then we 
have  $P\vdash \phi \in T(X)$.  We 
may form the element in $T(X)$
$$P_{\phi}:=(P\vdash \phi )\cdot \lambda^{-1};$$
this is a probability distribution. For by Theorem \ref{fourx}, its total is calculated as 
$\lambda ^{-1}$ multiplied  on 
$$\int _{X}1_{X}\; d(P\vdash \phi )(x) = \int _{X}1_{X}\cdot \phi 
(x)\; dP(x)=\int _{X}\phi (x)\; dP(x) = \lambda.$$ So we get $1$.

 Let us think of $\phi$ in the above consideration as a 
(generalized) ``event'' $A$, writing $A$ for $\phi$; also, let us 
write $P(B)$ for $\int_{X}B(x)\; dP(x)$, for general $B\in X\p R$. Then we have
$\lambda =P(A)$, and
the value of $P\vdash A$ on the ``event'' $B$ is $P(A\cdot B)$. 
Now $A\cdot B$ is the event $A\cap B$ (for the case of characteristic 
functions of subsets of $X$). So $P_{A}$ is $P(A\cap 
B)/P(A)$, the classical ``conditional probability of $B$ given $A$''.

\noindent Aarhus March 2011

\begin{verbatim}kock@imf.au.dk
\end{verbatim}


\begin{thebibliography}{99}

\bibitem{BET} R.\ Blute, T.\ Ehrhard and C.\ Tasson, A convenient 
differential category, to appear in Cahiers de Top.\ et G\'{e}om.\ 
Diff.\ Cat.

\bibitem{Bunge} M.\ Bunge, Relative functor categories and categories 
of algebras, {\em J.\ Algebra} 11 (1969), 64-101.

\bibitem{Borceux} F.\ Borceux, Handbook of Categorical Algebra 2, 
Cambridge University Press 1994.

\bibitem{Cramer}H.\ Cram\'{e}r, Sannolikhetskalkylen, Almqvist \& 
Wiksell Stockholm 1949.

\bibitem{deSchipper2}W.J.\ De Schipper, Commutative monads on 
symmetric closed categories, report nr.\ 34, Wiskundig Seminarium, Vrije 
Universiteit  Amsterdam 1974.

\bibitem{deSchipper1}W.J.\ De Schipper, Symmetric Closed Categories, 
Mathematical Centre Tracts 64, Amsterdam 1975.

\bibitem{EK} S.\ Eilenberg and  M.\ Kelly, Closed Categories, Proc.\ 
Conf.\ Categorical Algebra La Jolla 1965, 421-562, Springer Verlag 1966.

\bibitem{StoneSpaces} P.T.\ Johnstone,  Stone Spaces, Cambridge 
Studies in Advanced Mathematics 3, Cambridge University Press 1982.

\bibitem{Kelly} M.\ Kelly, Basic Concepts of Enriched Category 
Theory, London Math.\ Soc.\ Lecture Notes 64, Cambridge University 
Press 1982.

\bibitem{MSMCC}
A.\ Kock, Monads on symmetric monoidal closed categories.
  {\em Arch. Math. (Basel)}, 21:1--10, 1970.

\bibitem{DD}
A.\ Kock, On double dualization monads.
  {\em Math. Scand.}, 27:151--165, 1970.

\bibitem{CCGBCM} A.\ Kock, Closed categories generated by commutative monads.
  {\em J. Austral. Math. Soc.}, 12:405--424, 1971.

\bibitem{BCCM} A.\ Kock, Bilinearity and {C}artesian closed monads.
  {\em Math. Scand.}, 29:161--174, 1971.

\bibitem{SFMM} A.\ Kock, Strong functors and monoidal monads.
  {\em Arch. Math. (Basel)}, 23:113--120, 1972.

\bibitem{SDG} A.\ Kock, Synthetic Differential Geometry, London Math.\ Soc.\ 
Lecture Notes 51, Cambridge University 
Press 1981; Second Edition London Math.\ Soc.\ Lecture Notes 333, Cambridge University 
Press 2006.

\bibitem{PRSFA} A.\ Kock, Some problems and results in synthetic 
functional analysis, Category Theoretic Methods in Geometry Aarhus 
1983, Aarhus Var.\ Publ.\ Series 35 168-191, 1983.

\bibitem{MSPQ} A.\ Kock, Mathematical Structure of Physical Quantities, 
{\em Archive for Rational Mechanics and Analysis}, 107:99--104, 1989.

\bibitem{KM} A.\ Kriegl and P.\  Michor, The Convenient Setting of 
Global Analysis, 
Am.\ Math.\ Soc.\ 1997.

\bibitem{ACFPE} F.W.\ Lawvere, Algebraic Concepts on the Foundations 
of Physics and Engineering, MTH 461/561, Buffalo Jan.\ 1987.

\bibitem{CSQ} F.W.\ Lawvere, Categories of space and of quantity, 
 in: J. Echeverria et al (eds.), 
The Space of mathematics , de Gruyter, Berlin, New York (1992)

\bibitem{Lindner}H.\  Lindner, Affine parts of monads, {\em Arch.  
Math.\  (Basel)}, 33 (1979), 437-443.

\bibitem{Linton}F.\ Linton, Coequalizers in categories of algebras, 
Seminar on Triples and Categorical Homology Theory: 75-90, Springer 
Lecture Notes in Math. 80 (1969).

\bibitem{ManesBook}E.\ Manes, Algebraic Theories, Springer Graduate 
Texts in Math. 26, 1976.

\bibitem{MR} I.\ Moerdijk and G.E.\ Reyes, Models for Smooth 
Infinitesimal Analysis, Springer Verlag 1991.

\bibitem{Schwartz} L.\ Schwartz, M\'{e}thodes math\'{e}matiques pour les 
sciences physiques, Hermann Paris 1961.

\bibitem{Wraith} G.C.\ Wraith, Algebraic Theories (revised version), 
Aarhus Matematisk Institut Lecture Notes Series 22, 1975.

\end{thebibliography}
\end{document}